\documentclass[preprint,3pt]{elsarticle}

\usepackage{graphicx,amsmath,amsfonts}

\usepackage{times,amsfonts,amssymb,amsmath}
\usepackage{psfrag, color}
\usepackage{multirow}
\usepackage{comment}
\usepackage{caption}

\usepackage{color}

\def\c{{\bf c}}
\def\d{{\bf d}}

\def\f{{\bf f}}

\def\h{{\bf h}}
\def\i{{\bf i}}

\def\n{{\bf n}}

\def\p{{\bf p}}

\def\r{{\bf r}}

\def\t{{\bf t}}

\def\w{{\bf w}}

\def\z{{\bf z}}

\def\P{{\bf P}}

\def\U{{\bf U}}
\def\V{{\bf V}}
\def\K{{\bf K}}

\def\sign{{\rm sign}\;\!}

\def\Re{{\rm Re}\;\!}
\def\Im{{\rm Im}\;\!}
\def\CHI{{\bfm{\chi}}}
\def\signC{\zeta}

\newcommand{\be}{\begin{equation}}
\newcommand{\ee}{\end{equation}}
\newcommand{\ba}{\begin{eqnarray}}
\newcommand{\ea}{\end{eqnarray}}
\newcommand{\bi}{\begin{itemize}}
\newcommand{\ei}{\end{itemize}}

\newtheorem{rmk}{Remark}
\newtheorem{lma}{Lemma}

\newtheorem{thm}{Theorem}

\newtheorem{exm}{Example}[section]

\def\bfm#1{\boldsymbol{#1}}

\def\RR{\mathbb{R}}
\def\CC{\mathbb{C}}

\newcommand{\nrm}[1]{\left\| #1 \right\|}
\newcommand{\normC}[1]{\left| #1 \right|}
\newcommand{\dr}[1]{#1'}
\newcommand{\ddr}[1]{#1''}

\newcommand{\qqS}[1]{\mathcal{#1}^{\raise0.02ex\hbox{\scriptsize *}}}

\newcommand{\bern}[2]{B_{#2}^{#1}}

\newcommand{\abs}[1]{\left|#1\right|}
\newcommand{\conj}[1]{\overline{#1}}

\begin{document}

\begin{frontmatter}
\title{Construction of $G^2$ planar Hermite interpolants with prescribed arc lengths}


\author[UL,UR]{Marjeta Knez}
\address[UL]{FMF, University of Ljubljana, Jadranska 19, 1000 Ljubljana, Slovenia}
\address[UR]{IMFM, Jadranska 19, 1000 Ljubljana, Slovenia}
\ead{marjetka.knez@fmf.uni-lj.si}
\author[URome]{Francesca~Pelosi}
\address[Rome]{Dipartimento di Matematica,  Universit\`a di Roma ``Tor Vergata'', \\
Via della Ricerca  Scientifica, 00133 Roma, Italy}
\ead{pelosi@mat.uniroma2.it}
\author[US]{Maria Lucia Sampoli}
\address[US]{Department of Information Engineering and Mathematics, University of Siena, Italy}
\ead{marialucia.sampoli@unisi.it}

\begin{abstract}
In this paper we address the problem of constructing $G^2$ planar Pythagorean--hodograph (PH) spline curves, that interpolate points, tangent directions and curvatures, and have prescribed arc-length. The interpolation scheme is completely local. Each spline segment is defined as a PH biarc curve of degree $7$, which results in having a closed form solution of the $G^2$ interpolation equations depending on four free parameters. By fixing two of them to zero, it is proven that the length constraint can be satisfied for any data and any chosen ratio between the two boundary tangents. Length interpolation equation reduces to one algebraic equation with four solutions in general.
 To select the best one,  the value of the bending energy is observed.
Several numerical examples are provided to illustrate the obtained theoretical results and to numerically confirm that the approximation order is $5$.
\end{abstract}


\end{frontmatter}

\section{Introduction}
One of the fundamental problems of computer-aided design is to pass a curve through a given sequence of points. For example, in the design of cars and ships it is a standard practice to construct surfaces by first constructing networks of curves.
In addition to satisfy functional or aesthetic criteria, designed objects often have to exactly match not only a series of points but also derivatives. The problem can be efficiently addressed by $C^1$ or $C^2$ Hermite interpolation.
In industrial design, fair parametric curves and surfaces are the most preferred  representation to meet the requirements of design and modeling.
To this aim, it is sometimes preferable to consider a geometric continuity, that is simply the continuity of unit tangents
($G^1$ continuity), or even of the curvature function ($G^2$ continuity), see for instance \cite[Chapter~8]{FarinHoschekKim-02-Handbook}.
This kind of continuity is advantageous for many applications.
For instance, in NC milling applications, using spline curves with $G^2$ continuity produces motions with continuous accelerations, hence continuous cutting forces. This increases the possible speed of the manufacturing process, and/or the lifetime of the machinery.

On the other hand the possibility to construct curved paths satisfying given boundary conditions and with prescribed arc lengths is a fundamental problem in geometric design. Indeed such problems may arise in robot path planning,  computer animation, path planning for unmanned or autonomous vehicles, and related applications. In this context, it is well known that polynomial Pythagorean--hodograph curves possess polynomial cumulative arc length functions, so they are perfectly suited for the construction of exact solutions to such problems. {For instance the efficient use of PH curves in real-time control of CNC machines can be found in \cite{conway2012}, where PH curves are used for contour error computation of general free-form curved paths, while applications to robot motion path planning are presented in \cite{singh2018, su2018}.}

For planar PH curves, one of the first interpolation methods was given in \cite{FaroukiNeff'95}, where the interpolation of first order Hermite data was analysed. Later the problem was revisited in \cite{MeekWalton'97}. In \cite{ByrtusBastl'10} a characterization of the set of $G^1$ Hermite data for which planar PH cubic interpolants exist was given. The method of \cite{FaroukiNeff'95} was extended in \cite{Juttler'2000} to accommodate interpolation of $C^1$ data together with end–point curvatures, using planar PH curves of degree $7$, and it was shown that up to eight distinct interpolants exist. For quintic planar PH curves, several results on  first and  second order continuous spline interpolation can be found in the literature, see for instance \cite{FaroukiManniSest'03,FaroukiGiannSest'16, PSFM'07, Jaklic_etal'10, FaroukiHormNudo'20}.
In \cite{Jaklic_etal'14} a local $G^2$ interpolation with PH quintics that interpolate end points together with unit tangents and curvature values is considered. Although multiple solutions exist, a method to obtain a “good” interpolant is suggested through an asymptotic analysis.  However, the problem is highly non-linear and any general existence results seems to be unachievable.

The construction of curved paths with prescribed arc lengths satisfying given boundary conditions has received relatively little attention in the past, but it has recently been exposed as an important problem. The imposition of arc length constraints in the construction of PH curves was first considered in \cite{Huard_etal'14}, where a numerical scheme was employed to solve the system of non–linear equations that define a spatial $C^2$ PH quintic spline interpolating a sequence of nodal points with specified internodal arc lengths.
A closed–form solution to the problem of interpolating planar $G^1$ Hermite data under arc length constraints was developed in \cite{Farouki'16}, using planar PH quintics, and in \cite{Farouki'19} this approach was generalized to the spatial case. In \cite{Krajnc'17} rational PH curves are used to provide a closed-form solution to the problem of interpolating spatial $G^1$ data with prescribed arc lengths.

In this work we extend the interpolation problem analysed in \cite{Juttler'2000}, by adding the length
interpolation condition, while relaxing the $C^1$ to $G^1$ continuity conditions. In more detail, we consider the interpolation of two planar data points, two tangent directions and two curvatures by degree $7$ planar PH curves, that in addition have prescribed length.  The problem reduces to three non-linear equations with one free shape parameter. However, the analysis of existence of solutions is still too involved, and there exist data for which no interpolants can be found.  To make the problem theoretically as well as numerically simpler, we propose to replace a single PH curve with a PH biarc curve. Moreover, in order to reduce the number of free parameters coming from the splitting, we assume  that the biarc is $C^3$ continuous at the joint point. In this way we derive a closed form expression for the $G^2$ interpolant, which still depends on four free parameters. Then we fix two of them to zero by requiring that the first and the second order derivatives at  boundary points are orthogonal, and set the other two proportional to the lengths of boundary tangents. Eventually we are able to prove that the length interpolation condition, that reduces to one simple algebraic equation, has a solution for any data and any value of the remaining shape parameter. Beside of being easy to implement and simple to use in practise, the proposed method  can be directly applied to a (local) construction of $G^2$ continuous interpolating splines. {To measure the fairness of the resulting  interpolants, we observe the value of the curve bending energy. To visualize the results, we use the porcupine plot of the curvature (see for instance \cite{Fairshape}).  The porcupine curvature plot consists of lines emanating from points on the curve in the direction of $\sign{(\kappa)} \, \n$ (where $\n$ is a unit normal) with lengths proportional to the magnitude of the curvature $\kappa$. In this way it is possible also to see immediately the curvature continuity and its sign changes.
}

The remainder of this paper is organized as follows. In the next section we introduce the notation and give some basic definitions. Then in Section~\ref{sec:probl} we present the problem and propose a possible approach to solve it. We show in particular that a solution could not exist for some given data with a single PH curve of degree $7$.
Section~\ref{sec:biarc} presents a  biarc construction and its theoretical analysis.
In Section~\ref{sec:Examples} some numerical examples are given to illustrate the performance of the proposed method. Finally, Section ~\ref{sec:conclusion} briefly summarizes the contributions of the present study, and identifies possible directions for further investigation.

\section{Preliminaries}
The distinctive property of a polynomial planar PH curve ${\r}(t)=\left(r_1(t),r_2(t)\right)$ is that its hodograph
${\r}'(t)=(\dr{r_1}(t),\dr{r_2}(t))$ satisfies the Pythagorean condition
\be
\label{PH_cond}
\dr{r_1}^{\,2}(t) \,+\, \dr{r_2}^{\,2}(t) \,=\, \sigma^2(t)
\ee
for some polynomial $\sigma(t)$, which specifies the parametric speed of $\r(t)$, i.e. the derivative of the arc length $s$ with respect to the curve parameter $t$. This feature provides planar PH curves with many attractive computational properties: they have rational unit tangents and normals, curvatures, and offset curves. Moreover their arc lengths are exactly computable so that they are ideally suited to real-time precision motion control applications, \cite{Farouki-PH-book-2008}. In the complex representation, \cite{Farouki'94},  a planar PH curve segment is constructed from a complex polynomial ${\w}(t)= u(t)+ \i \,v(t)$ by integrating the expression
\be\label{rprime}
\r'(t)={\w}^2(t).
\ee
We use bold characters to interchangeably denote planar points/vectors and complex
numbers, while real (scalar) quantities are denoted by italic characters.  {Note also that we identify complex numbers with points/vectors in $\RR^2$. Thus, the curve $\r(t)$ is further expressed as $\r(t) = r_1(t) + \i\, r_2(t)$, and similarly, all the interpolation points/vectors will be given as complex numbers.}

If $\w(t)$ is of degree $m$, it may be defined by its complex Bernstein coefficients $\w_0,\ldots,\w_m$,
\be
\w(t)=\sum_{i=0}^m \w_i B_i^m(t), \;\; B_i^m(t)={{m}\choose{i}} (1-t)^{m-i}\ t^i.
\ee
The above expression is also called {\it  complex preimage curve} of $\r$.
Integrating (\ref{rprime}) yields a planar PH curve of degree $n = 2m + 1$. Denoting the hodograph (\ref{rprime}) with $\h(t)$, it can be written in  B\'ezier form as
\be \label{h}
\h(t)=\sum_{i=0}^{2m} \h_i B_i^{2m}(t),
\ee
where the coefficients are related to the complex polynomial $\w(t)$ as
\be \label{h_expr}
\h_i=\sum_{j=\max(0,i-m)}^{\min(m,i)}
\frac{{{m}\choose{j}} {{m}\choose{i-j}}} {{{2m}\choose{i}}}
 \w_j \ \w_{i-j}, \;\;i=0,\ldots,2m.
\ee
The planar PH curve of degree $n$  can be also expressed in the B\'ezier representation
\be \label{PH}
\r(t)=\sum_{i=0}^{n} \p_i B_i^{n}(t),
\ee
and taking into account that $\h(t)=\r'(t)$  we have
\be
\h(t)=n \ \sum_{i=0}^{n-1} \Delta \p_i B_i^{n-1}(t), \quad  \Delta \p_i:= \p_{i+1}-\p_i,
\ee
and
\be \label{ph_cp}
\p_i=\p_{i-1}+\frac{1}{n}\h_{i-1},\; i=1,\ldots,n,
\ee
where $\p_0$ is a free integration constant. From (\ref{rprime}) the parametric speed $\sigma(t)$, unit tangent $\t(t)$ and curvature $\kappa(t)$ of the curve $\r(t)$ may be formulated in terms of $\w(t)$ as, see \cite{Farouki'94},
\be
\sigma(t)=| \w(t)|^2 ,\; \quad \t(t)=\frac{\w(t)^2}{\sigma(t)}, \; \quad \kappa(t)=2\ \frac{\Im(\overline{\w}(t)\ \w'(t))}{\sigma^2(t)},
\ee
where $|\cdot|$ denotes the absolute value of a complex number, and $\overline{\w}$ is a conjugate of $\w$.
Since $\sigma(t)$ is a polynomial of degree $2m$, the cumulative arc length function
$$
s(t)=\int_0^t \sigma(\lambda)\ d\lambda$$
is likewise just a polynomial of degree $2m+1$.
\begin{rmk} \label{rem-normals}
A normal vector at each point on the curve $\r$ is perpendicular to the unit tangent $\t$. It can be computed  by a $90^{\circ}$ clockwise or counterclockwise rotation of $\t$. In what follows we choose the
later one, i.e. $\n = {\t}^{\perp}$, where
$(x,y)^{\perp}:= (-y,x)$.
\end{rmk}

The quintics, obtained by integrating ${\bf r}'(t)={\bf w}^2(t)$ with a complex quadratic ${\bf w}(t)$, are the lowest--order PH curves with sufficient shape flexibility for free--form design.
However, for our purposes they do not possess enough degrees of freedom to fulfil the $G^2$ interpolation at the boundary, and the arc-length interpolation. Hence in this paper we shall focus on PH curves of degree 7.\\
In this case, by choosing a cubic polynomial $\w(t)$ with Bernstein coefficients $\w_0,\ \w_1,\ \w_2, \ \w_3$, the expression (\ref{h_expr}) reduces to
\begin{equation}\label{h6}
\begin{split}
\h_0 &= \w_0^2, \quad 
\h_1= \w_0 \w_1, \quad 
\h_2= \frac{1}{5}(2\w_0\w_2+3\w_1^2), \quad 
\h_3=\frac{1}{10}(\w_0\w_3+9\w_1\w_2), \\
\h_4&=\frac{1}{5}(2\w_1\w_3+3\w_2^2), \quad 
\h_5=\w_2 \w_3, \quad 
\h_6= \w_3^2. 
\end{split}
\end{equation}
The parametric speed is a polynomial of degree $6$, which can be expressed in  the B\'ezier form as
\be
\sigma(t)= \sum_{i=0}^6 \sigma_i B_i^6(t)
\ee
with coefficients
\begin{equation}\label{sigma}
\begin{split}
\sigma_0 &= \abs{\w_0}^2, \quad 
\sigma_1 = \Re\left(\w_0\  \conj{\w}_1\right), \quad 
\sigma_2 = \frac{1}{5}\left(2 \Re(\w_0 \ \conj{\w}_2) + 3 \abs{\w_1}^2\right), \\ 
\sigma_3 & =\frac{1}{10}\Re\left(\w_0\ \conj{\w}_3 + 9 \w_1\ \conj{\w}_2\right), \\
\sigma_4&=\frac{1}{5}\left(2 \Re(\w_1 \ \conj{\w}_3) + 3 \abs{\w_2}^2\right), \quad
\sigma_5=\Re\left(\w_2\  \conj{\w}_3\right), \quad 
\sigma_6= \abs{\w_3}^2. 
\end{split}
\end{equation}

The total arc length is
\be \label{L}
 L_{\r} =
 s(1)=\frac{\sigma_0+\sigma_1+ \ldots +\sigma_6}{7}.
\ee
\section{The interpolation problem}\label{sec:probl}
Given two data points $\P_0$ and $\P_1$, two associated tangent directions $\t_0$ and $\t_1$ (with $\normC{\t_0}=
\normC{\t_1}=1$) and two signed curvatures
$\kappa_0$, $\kappa_1$, the problem we want to address is to find a PH curve $\r: [0,1] \rightarrow \CC$ of degree 7 that interpolates the data
\begin{equation} \label{eq:intCond}
\r(i) = \P_ i, \quad \frac{\dr{\r}(i)}{|\dr{\r}(i)|} = \t_i \, \quad \kappa(i) =\kappa_i, \quad i=0,1,
\end{equation}
and has the length equal to a given $L$, such that $L>|{\Delta \P_0}|$.
It is easy to verify that these $G^2$ interpolation conditions are satisfied if the control points are equal to (see \cite{Lu'15})
\begin{equation} \label{G^2_cond}
\begin{split}
\p_0 &= \P_0, \quad \p_1 = \P_0+\frac{\alpha_0^2}{7} \t_0, \quad \p_2 = \P_0 + \left(\frac{2 \alpha_0^2}{7}+\frac{\beta_0}{42}\right) \t_0+\kappa_0 \frac{\alpha_0^4}{42} \n_0, \\
\p_7 &= \P_1, \quad \p_6 = \P_1 - \frac{\alpha_1^2}{7} \t_1,\quad \p_5 = \P_1 + \left(-\frac{2 \alpha_1^2}{7}+\frac{\beta_1}{42}\right) \t_1 + \kappa_1 \frac{\alpha_1^4}{42} \n_1,
\end{split}
\end{equation}
for some  $\alpha_0, \alpha_1$ and $\beta_0,\beta_1$, where $\n_0=\i\, \t_0$, $\n_1=\i\, \t_1$ are the unit normals at the end-points (see Remark~\ref{rem-normals}).
Of course, these parameters need to be computed so that $\r$ satisfies the PH condition.

From the expressions for $\p_1$ and $\p_6$ of \eqref{G^2_cond} and relations
\eqref{ph_cp} and \eqref{h6} 
we get the following complex quadratic equations
\be \label{end_tang_cond}
\w_0^2=\alpha_0^2\ \t_0, \;\quad \w_3^2=\alpha_1^2\ \t_1.
\ee
The next lemma follows by the elementary computations.
\begin{lma} \label{lemma1}
Let $\c:=c_1+ \i\, c_2 \in \mathbb{C}$, with $\c \neq0$. Then the equation $\z^2=\c$ has two solutions
$$
\z= \signC \,\CHI(\c), \quad \CHI(\c) := \frac{\sqrt{2}}{2}\left(\sqrt{\abs{\c}+c_1} + \i \frac{c_2}{\sqrt{\abs{\c}+c_1}}\right)
,\quad  \signC\in\{-1,1\}.
$$
\end{lma}
Using Lemma~\ref{lemma1} we can solve the equations (\ref{end_tang_cond}) in terms of the preimage coefficients:
\be \label{eq1}
\w_0= \alpha_0\ \zeta_0 \  \CHI(\t_0) \quad \quad \w_3= \alpha_1\ \zeta_1 \ \CHI(\t_1), \quad \zeta_0,\zeta_1\in \{-1,1\}.
\ee
Without loosing generality we set $\zeta_0= \zeta_1=1$ and observe positive as well as negative values  for $\alpha_0$ and $\alpha_1$. Now, considering the expressions for $\p_2$ and $\p_5$ of \eqref{G^2_cond}, the relations in \eqref{h6} for $\h_1$ and $\h_5$, and \eqref{ph_cp} we obtain
\begin{equation} \label{eq2}
\w_1 =\frac{1}{\w_0} \left( \left(\alpha_0^2+\frac{\beta_0}{6}\right) \t_0 + \kappa_0 \frac{\alpha_0^4}{6} \n_0 \right), \quad \w_2 =\frac{1}{\w_3}\left( \left(\alpha_1^2 - \frac{\beta_1}{6}\right) \t_1 - \kappa_1 \frac{\alpha_1^4}{6} \n_1 \right).
\end{equation}
From  \eqref{eq1} and \eqref{eq2} we see that using the end tangents and end curvature conditions we have been able to express all the coefficients of the preimage function $\w(t)$ in terms of the free parameters $\alpha_0, \alpha_1,$ $\beta_0$ and $\beta_1$.
To complete the interpolation conditions (\ref{eq:intCond}) we have to impose the end point interpolation:
$$
7\ \Delta \P_0 = \sum_{i=0}^6 \h_i, \quad  \; \Delta \P_0= \P_1-\P_0.
$$
From (\ref{h6}), this results in the following quadratic equation
\begin{align} \label{eq3}
10(\w_0^2 + \w_3^2) +  6(\w_1^2 +  \w_2^2) &+ 10(\w_0\ \w_1 + \w_2 \ \w_3)+
 \w_0\ \w_3 + \\ \nonumber
 & 4(\w_0\ \w_2 + \w_1\ \w_3)+ 9\w_1\ \w_2 = 70 \Delta \P_0.
\end{align}
At last, from (\ref{sigma}) and  (\ref{L}), requiring the interpolant to have a specified  length $L$ yields the condition
\begin{align}\label{eq4}
10(| \w_0 |^2 + | \w_3 |^2) +  6(| \w_1|^2 + | \w_2|^2)+ &\Re\left( 10(\w_0\  \conj{\w}_1+ \w_2\  \conj{\w}_3)+
\w_0\ \conj{\w}_3 \right.\\ \nonumber
& \left. 4(\w_0 \ \conj{\w}_2+\w_1 \ \conj{\w}_3)
+9 \w_1\ \conj{\w}_2\right)   = 70 L.
\end{align}
Equations \eqref{eq3} and \eqref{eq4} give three (highly) non-linear scalar equations for four unknowns
$\alpha_0$, $\alpha_1$, $\beta_0$ and $\beta_1$, thus the solutions depend on one extra parameter. This one  degree of freedom can be fixed by prescribing the ratio between lengths of boundary tangents.

We note that if we omit the length constraint \eqref{eq4} and fix the parameters $\alpha_0$ and $\alpha_1$ by assigning the lengths of both boundary tangent vectors, we obtain a different interpolation problem, which has been studied by J\"uttler in  \cite{Juttler'2000},  resulting with $C^1/G^2$ PH interpolants. More precisely, it is shown in \cite{Juttler'2000} that
for fixed $\alpha_0$ and $\alpha_1$, there exist  up to eight possible solutions of equations \eqref{eq3} that come as roots of two bivariate quartic polynomials.  Thus, any general result on the existence of solutions seems to be impossible, so a comprehensive asymptotic analysis is provided instead and it is shown that for data taken from a smooth parametric curve, defined on some small interval, there exist solutions having the approximation order $6$.

Relaxing the $C^1$ boundary conditions to $G^1$, as done in \eqref{eq:intCond}, and adding the equation \eqref{eq4} makes the interpolation problem even more non-linear and more difficult to analyse. In addition, it is easy to come up with examples for which no solutions exist or the solutions are not visually pleasing, as demonstrated in the next example. {To numerically solve the nonlinear equations in this example we have used the program package {\tt Mathematica} and its function {\tt Solve}, that gives all the solutions (including the complex ones) of the polynomial system of equations.}
\begin{exm}  \label{exm-OnePHCurve}
Let the interpolation data be chosen as
\begin{equation} \label{example-int-data}
\bfm{P}_0 = 0+ 0 \, \i,\quad  \bfm{P}_1 = 1 + 0\,  \i, \quad \bfm{t}_0 = \cos\theta_0 + \i \sin \theta_0,\quad
\bfm{t}_1 = \cos\theta_1 + \i  \sin \theta_1
\end{equation}
for $\theta_0 = -\frac{\pi}{4}$ and $\theta_1 = -\frac{\pi}{8}$. Moreover, let
$\kappa_0 = 1$, $\kappa_1 = -1$, $L=1.1$, and let us additionally require that $|\r'(0)|=|\r'(1)|$. Then there  exist two PH interpolants shown in Fig.~\ref{fig:example-int-data} together with control polygons and porcupine curvature plots.
\begin{figure}[htb]
\centering
\begin{minipage}{0.45\textwidth}
\includegraphics[width=1\textwidth]{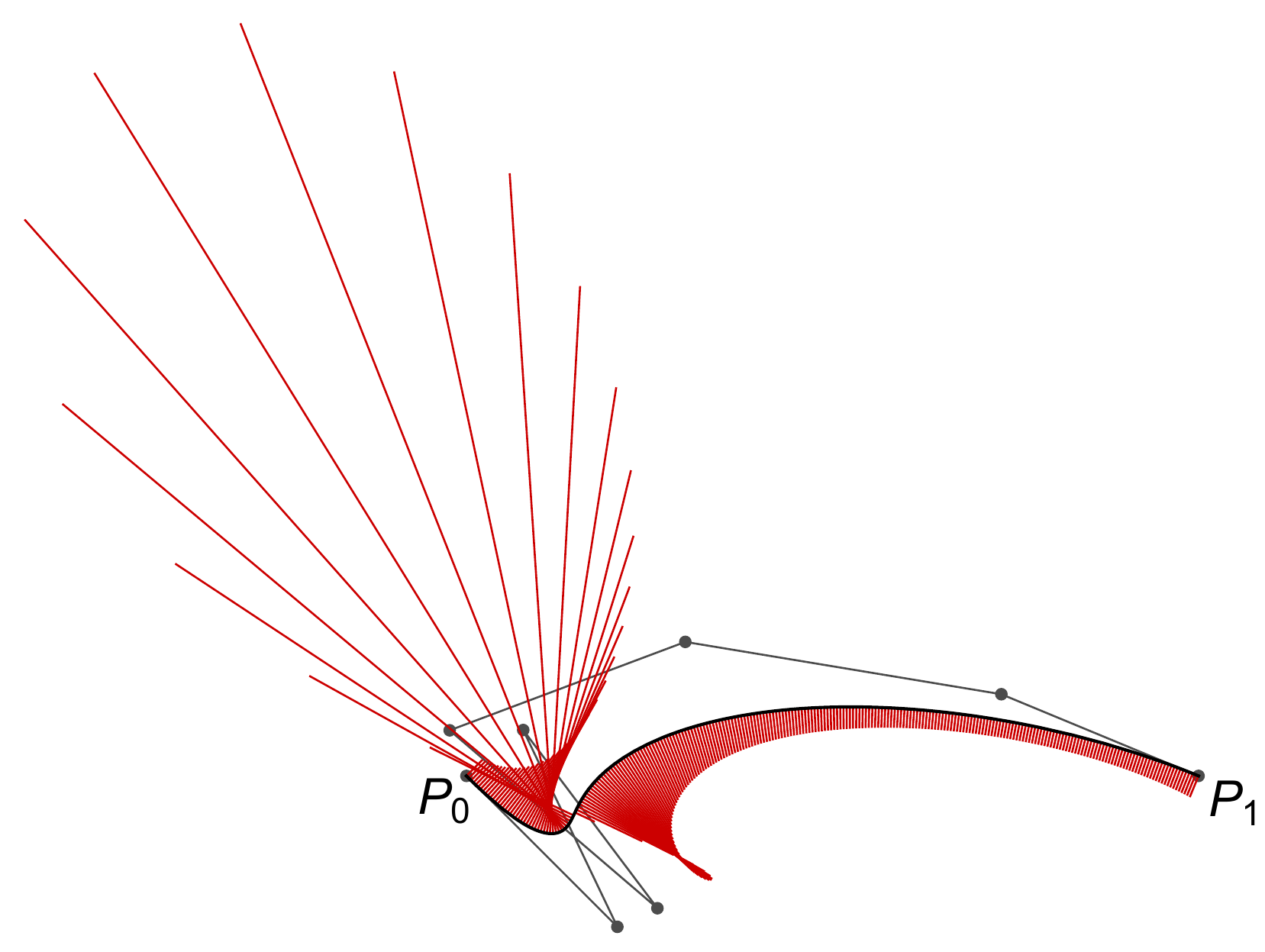}
\end{minipage}
\hskip1cm
\begin{minipage}{0.45\textwidth}
\includegraphics[width=0.9\textwidth]{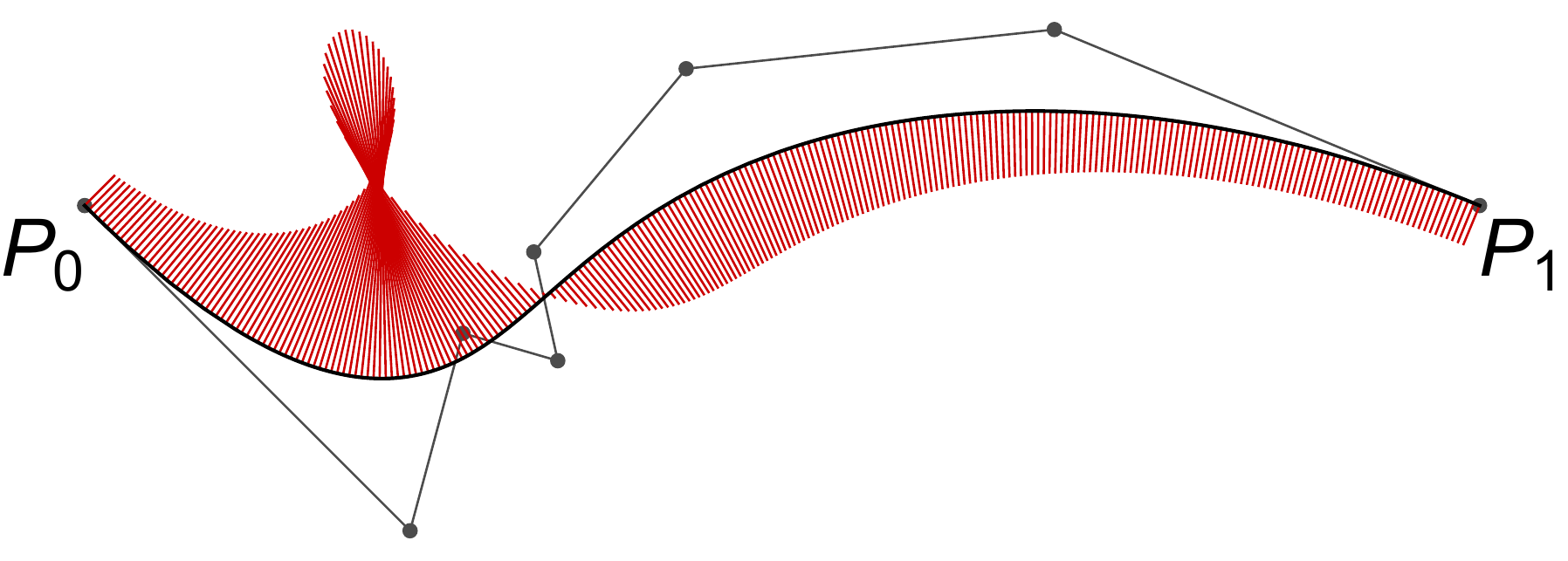}
\end{minipage}
\caption{Two PH interpolants of degree $7$ from Example \ref{exm-OnePHCurve}  together with the control polygons and porcupine curvature plots (with proportional factor equal to 0.03).}
\label{fig:example-int-data}
\end{figure}
If we decrease $L$ to $1.05$ (or less), no real solutions of \eqref{eq3} and \eqref{eq4} exist.
Furthermore, changing the sign of $\theta_1$ and $\kappa_1$, we get convex data for which there exist four different PH interpolants, but none of them has a nice shape, as shown in Fig.~\ref{fig:example-int-data-a}.
\begin{figure}[htb]
\centering
\begin{minipage}{0.24\textwidth}
\includegraphics[width=1\textwidth]{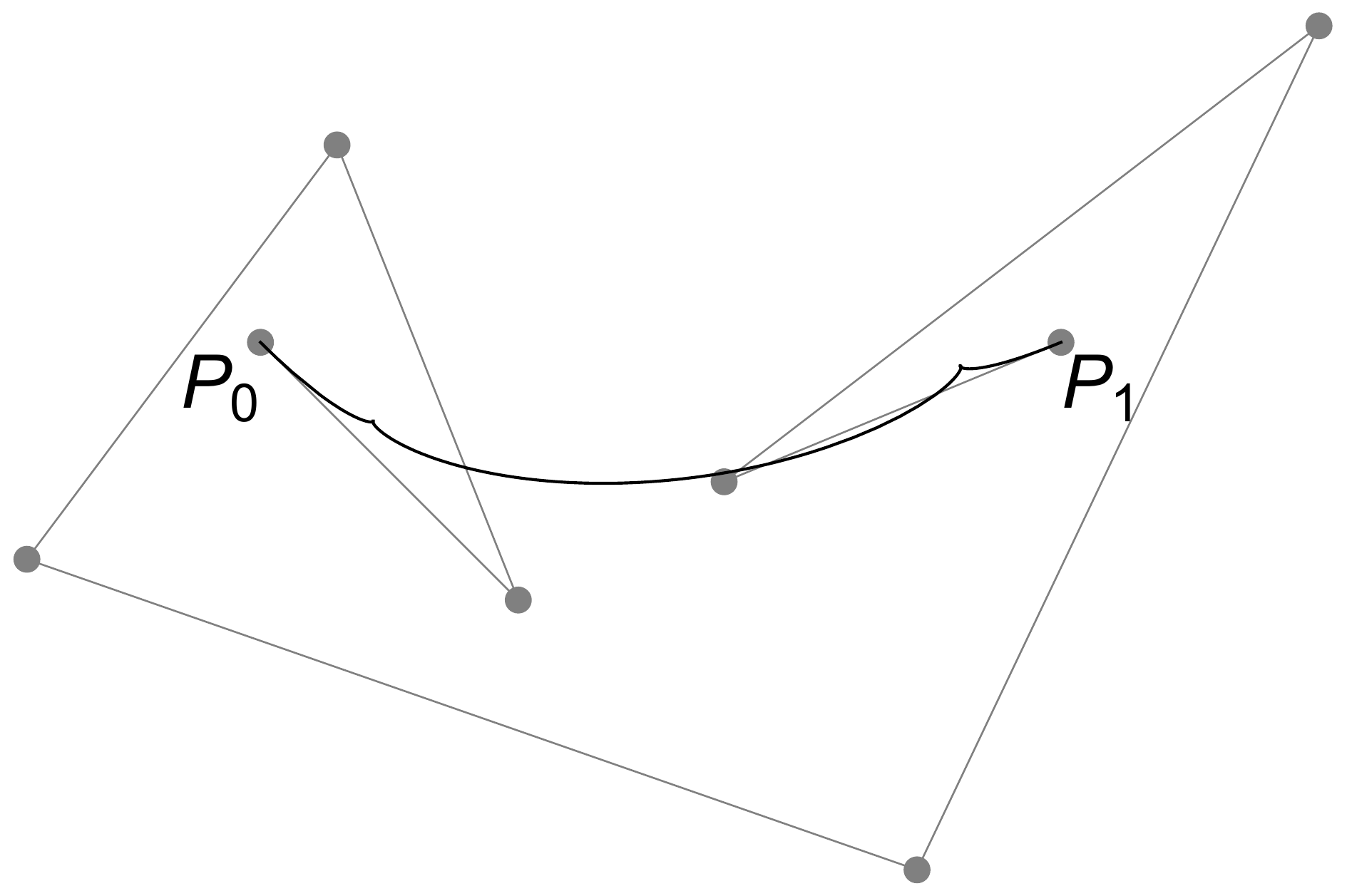}
\end{minipage}
\begin{minipage}{0.24\textwidth}
\includegraphics[width=1\textwidth]{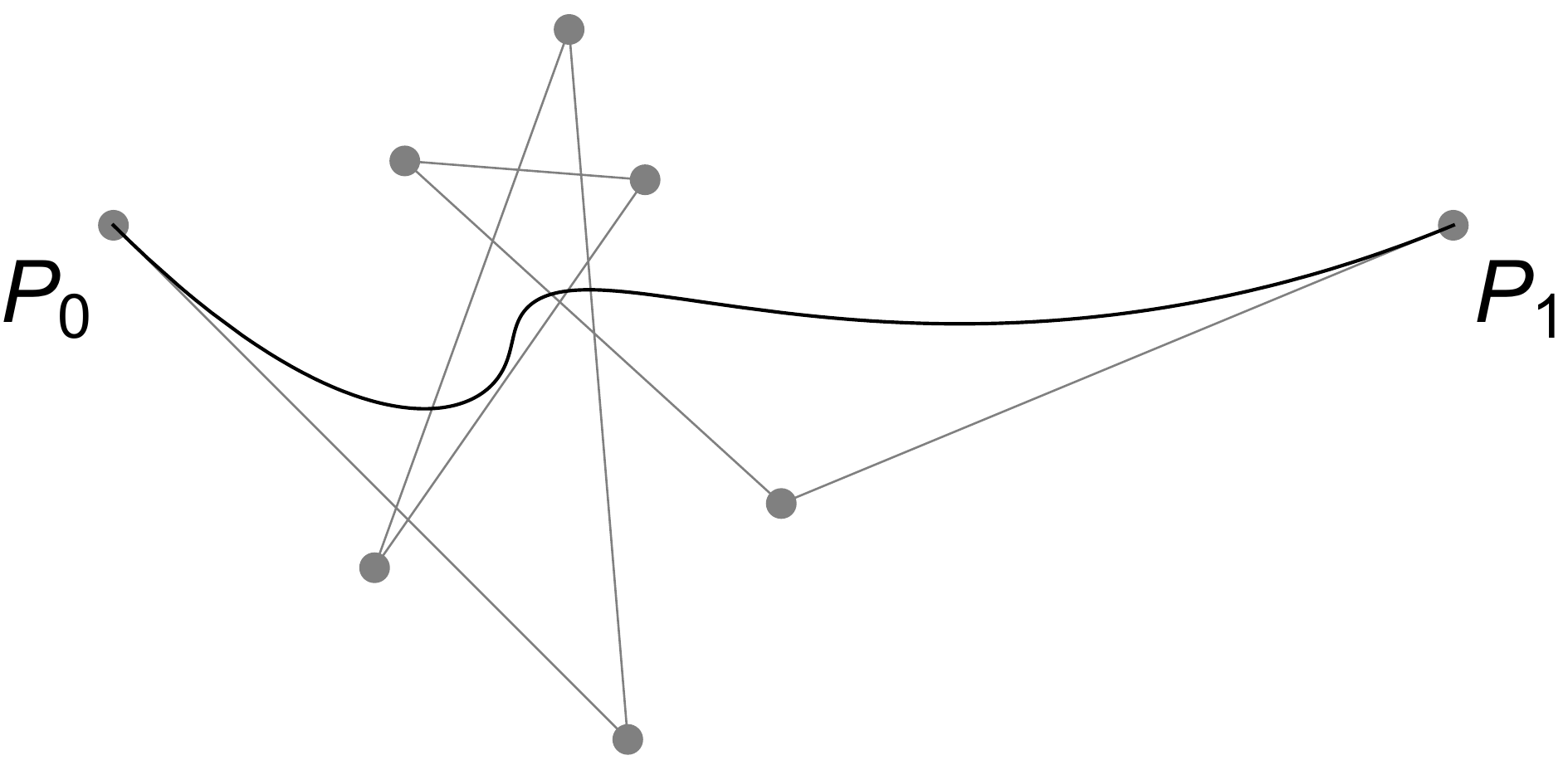}
\end{minipage}
\begin{minipage}{0.24\textwidth}
\includegraphics[width=1\textwidth]{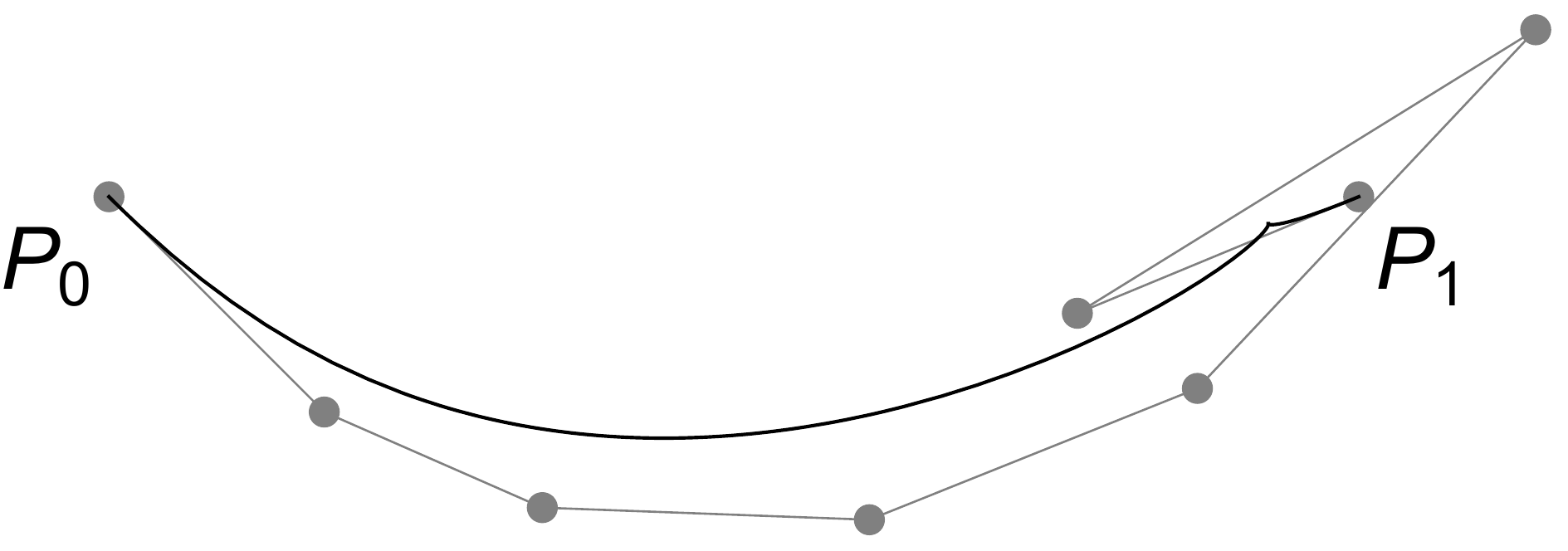}
\end{minipage}
\begin{minipage}{0.24\textwidth}
\includegraphics[width=1\textwidth]{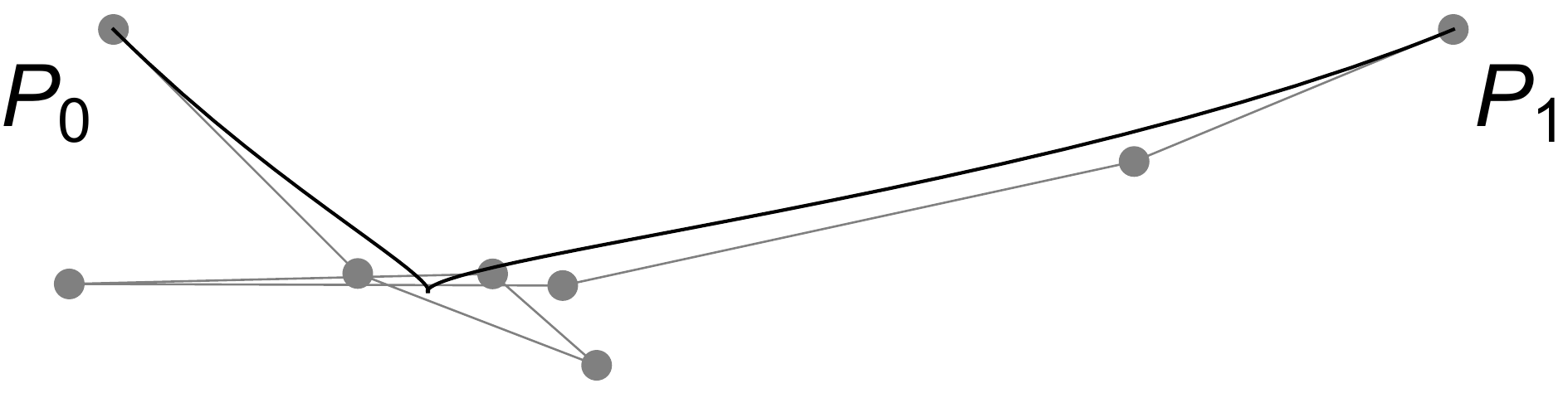}
\end{minipage}
\caption{Four PH interpolants $\r$ of degree $7$ for the data \eqref{example-int-data} with
$\theta_0 = -\frac{\pi}{4}$, $\theta_1 = \frac{\pi}{8}$, $\kappa_0 = \kappa_1 = 1$, $L=1.1$, and
 $|{\r'(0)}|=|{\r'(1)}|$.
}
\label{fig:example-int-data-a}
\end{figure}
\end{exm}

The main difficulty in analysing  equation \eqref{eq3} is that the unknowns are not the coefficients of the preimage curve, but they are parameters that come from $G^2$ conditions. One way to overcome this problem would be to raise the degree of the PH curve to $9$, which would provide us with one additional free preimage coefficient and much simpler solution of $G^2$ continuity equations. An alternative approach, examined in this paper,  is to leave the degree equal to $7$, but replace the polynomial curve with a biarc one.  As a consequence, the solution of the $G^2$ continuity equations can be derived in a closed form. Furthermore, by an appropriate choice of the additional free parameters, it can be proven that the length constraint \eqref{eq4} can always be fulfilled. Details are given in the following section.

\section{Biarc construction}\label{sec:biarc}
In order to address the Hermite interpolation constraints as well as to match the prescribed length for any arbitrary set of data, 
we consider the construction of a PH biarc segment.

A PH biarc of degree 7 can be defined as a piecewise curve
\begin{subequations} \label{def-biarc}
\begin{equation} \label{def-biarc-1}
\r(t) = \begin{cases}
\r_A(t), & t \in  \left[0, \frac{1}{2}\right)\\
\r_B(t), & t \in  \left[\frac{1}{2},1\right]
\end{cases},
\end{equation}
having each segment
expressed in the B\'ezier form as
\begin{equation}  \label{def-biarc-2}
\r_A(t) = \sum_{i=0}^{7}\p_{A,i} \bern{7}{i}(2t),  
\quad
\r_B(t) = \sum_{i=0}^{7}\p_{B,i} \bern{7}{i}(2t-1).   
\end{equation}
\end{subequations}
Considering their hodographs
\begin{equation}\label{eq-biarcHod}
\h_A(t)=\sum_{i=0}^{6} \h_{A,i} B_i^{6}(2t),\quad \; \h_B(t)=\sum_{i=0}^{6} \h_{B,i} B_i^{6}(2t-1),
\end{equation}
it holds that  $\h_A(t)=\w_A^2(t)$ and $\h_B(t)=\w_B^2(t)$, where the preimages are expressed in terms of cubic complex polynomials
\be
\w_{A}(t) = \displaystyle{\sum_{i=0}^{3}\w_{A,i} \bern{3}{i}(t)}, \; \quad
\w_{B}(t) = \displaystyle{\sum_{i=0}^{3}\w_{B,i} \bern{3}{i}(t)}.
\ee
Hence, from \eqref{h} and \eqref{h6}, the hodograph of the PH biarc curve of degree 7 is completely defined once the complex coefficients $\w_{A,i}$, $\w_{B,i}$, $i=0,1,2,3$,
are determined. Then, by integrating \eqref{eq-biarcHod}, we  obtain the control points of each PH segment
\eqref{def-biarc-2}. Choosing the free integration constants as $ \r_{A,0} = \P_0$ and $\r_{B,7} = \P_1$, we have
$$
\p_{A,i+1} = \P_0 + \frac{1}{14}\sum_{j=0}^{i} \h_{A,j},\;\;i=0, 1, \dots,  6,
\quad  \p_{B,i} = \P_1 - \frac{1}{14}\sum_{j=i}^{6} \h_{B,j},\;\; i=6,5, \dots, 0.
$$
Thus, the interpolation conditions $\r(0)=\P_0$ and $\r(1)=\P_1$ are already
achieved.
Let us now impose the $G^1$ conditions. Similarly to \eqref{eq1}, we can determine
\begin{subequations}\label{eq-W-biarc}
\begin{equation} \label{eq1_biarc}
\w_{A,0}= \alpha_0\, \signC_0\,  \CHI(\t_0),  \quad \w_{B,3}=  \alpha_1\, \signC_1 \, \CHI(\t_1)
\end{equation}
for four different sign choices $\left(\signC_0,\signC_1\right) \in \left\{(1,1), (1,-1), (-1,1), (-1,-1)\right\} $.
Again, we can fix
$\zeta_0 = \zeta_1 = 1$ and further observe positive as well as negative values of parameters $\alpha_0$ and $\alpha_1$.
Analogously to \eqref{eq2}, we obtain from the $G^2$ conditions that
\begin{align} \label{eq2_biarc}
\w_{A,1} &=\frac{1}{\w_{A,0}} \left( \left(\alpha_0^2+\frac{\beta_0}{12}\right) \t_0 + \kappa_0 \frac{\alpha_0^4}{12} \n_0 \right), \nonumber\\
\\
\w_{B,2} &=\frac{1}{\w_{B,3}}\left( \left(\alpha_1^2 - \frac{\beta_1}{12}\right) \t_1 - \kappa_1 \frac{\alpha_1^4}{12} \n_1 \right), \nonumber
\end{align}
where we should note that, in comparison to \eqref{eq2}, the different coefficients 
are due to the fact that the biarc segments are defined over halved intervals.

What is left is to assure that $\r$ is $G^2$ continuous at the joint parameter $t= \frac{1}{2}$.  As in the case of macro-elements, the usual approach to reduce the number of free parameters, coming from splitting, to the number of free parameters needed to get a simple interpolation construction, is to require additional smoothness at the joint point. In our case, requiring
$C^3$ continuity of $\w_A$ and $\w_B$ at $t=\frac{1}{2}$ would lead to a polynomial curve $\r$, i.e. the curve examined in the previous section. Thus, to obtain the additional freedom needed while keeping the interpolation scheme as simple as possible,  we require that the preimage biarc is $C^2$ continuous, which implies the biarc \eqref{def-biarc} to be $C^3$ continuous provided $\r_A(\frac{1}{2})=\r_B(\frac{1}{2})$ holds true.

Recalling the geometrical construction in terms of control points for achieving $C^2$ continuity of $\w_A$ and $\w_B$ at
$t=\frac{1}{2}$, we introduce a new control point, expressed with a complex number $\d\in \CC$, and define
%
\begin{equation} \label{eq-W-biarc-2}
\w_{A,2}=\frac{1}{2}(\w_{A,1}+\d), \quad \w_{B,1}=\frac{1}{2}(\w_{B,2}+\d)
\end{equation}
and
\begin{equation} \label{eq-W-biarc-3}
\w_{A,3}=\w_{B,0}=\frac{1}{2}(\w_{A,2}+\w_{B,1})=\frac{1}{4} \w_{A,1}+\frac{1}{2}\d+\frac{1}{4}\w_{B,2}.
\end{equation}
\end{subequations}
It remains to impose the $C^0$ continuity at the joint point, $\r_A(\frac{1}{2})=\r_B(\frac{1}{2})$. Hence we must choose $\d$ so that
\begin{equation} \label{eq-C0-cond}
\Delta\P_0=\frac{1}{14} \sum_{i=0}^{6} (\h_{A,i}+\h_{B,i}).
\end{equation}
Now considering the relations, given in \eqref{h6}, between the coefficients of the hodograph and those of the preimage for each subsegment,  after some computation, we can rewrite the equation \eqref{eq-C0-cond} in terms of $\d$, $\w_{A,0}$, $\w_{A,1}$, $\w_{B,2}$  and $\w_{B,3}$ as
\be \label{quad_eq}
\d^2+  2\ \U\ \d+ \V=0
\ee
with
\begin{align}
52\ \U &=   5\ (\w_{A,0}+\w_{B,3})+39\ (\w_{A,1}+\w_{B,2}), \nonumber\\
52\ \V &= 40\ (\w_{A,0}^2+\w_{B,3}^2) +49\ (\w_{A,0}\w_{A,1}+\w_{B,2}\w_{B,3}) +62 \ (\w_{A,1}^2+\w_{B,2}^2)  \label{def-of-UandV}\\
   & \hskip3.5cm   +\  \w_{A,0}\w_{B,2} +\w_{A,1}\w_{B,3}+ 28\w_{A,1}\w_{B,2} -560\ \Delta\P_0. \nonumber
\end{align}
 The quadratic complex equation \eqref{quad_eq} can be rewritten as $\left(\d+ \U \right)^2= \U^2 -\V$, and so using Lemma \ref{lemma1},
 we obtain an explicit solution for $\d$,
 \begin{equation}\label{eq-sol-d}
 \d= \signC_{\d}\, \CHI\left( \U^2 -\V \right)- \U, \quad \signC_{\d} \in \{-1,1\},
 \end{equation}
  in terms of $\w_{A,0}$, $\w_{A,1}$, $\w_{B,2}$  and $\w_{B,3}$, which depend on the free parameters $\alpha_0$, $\alpha_1$, $\beta_0$ and $\beta_1$ (see \eqref{eq1_biarc}, \eqref{eq2_biarc}).
The results are summarized in the next theorem.
\begin{thm} \label{thm-main-1}
A planar PH biarc \eqref{def-biarc} of degree $7$ that interpolates $G^2$ data \eqref{eq:intCond} is given in a closed form, that follows from the preimage, given by
\eqref{eq-W-biarc} and \eqref{eq-sol-d}, and depends on four free parameters $\alpha_0$, $\alpha_1$, $\beta_0$, $\beta_1$ and the sign choice $\signC_{\d} \in \{-1,1\}$.
\end{thm}


The last step is to require that the resulting biarc curve has the prescribed length $L$.
Taking into account \eqref{L}, in the biarc case, this constraint leads to the equation
\begin{equation*}
L= \frac{1}{14} \sum_{i=0}^{6} (\sigma_{A,i}+\sigma_{B,i}),
\end{equation*}
where $\sigma_{A,i}$ and $\sigma_{B,i}$ are the B\'ezier ordinates of $\sigma_A(t)=\w_A(t) \overline{\w}_A(t)$ and
 $\sigma_B(t)=\w_B(t) \overline{\w}_B(t)$, respectively.
 By using \eqref{sigma}, \eqref{eq-W-biarc} and \eqref{eq-sol-d}, we obtain, with some computation, the equation
 \begin{equation*}
   \abs{\d}^2 + 2 \, \Re\left(\d\, \conj{\U}\right) +v=0,
 \end{equation*}
 where
 \begin{equation}\label{v}
 \begin{split}
 52\ v = & 40 \abs{\w_{A,0}}^2 +  40 \abs{\w_{B,3}}^2 + 62 \abs{\w_{A,1}}^2 + 62 \abs{\w_{B,2}}^2 \\
 & +  \Re\left(49 \w_{A,0} \conj{\w}_{A,1} +49 \w_{B,2} \conj{\w}_{B,3}
 + 28 \w_{A,1} \conj{\w}_{B,2} + \w_{A,0} \conj{\w}_{B,2} +\w_{A,1} \conj{\w}_{B,3}
  \right) -560\ L.
 \end{split}
 \end{equation}
 If we denote $\K=\signC_{\d} \, \CHI\left(\U^2 -\V\right)$, then $\abs{\K}^2 = \abs{ \U^2 -\V}$
and
 \begin{align*}
 & \abs{\d}^2 +2\  \Re\left(\d\, \conj{\U}\right) =
 \left(\K - \U\right)
\left(\conj{\K} - \conj{\U}\right) +  2\ \Re\left(\K \conj{\U}\right) -2 \  \abs{\U}^2\\
& = \abs{\K}^2 - 2\ \Re\left(\K \conj{\U}\right) + \abs{\U}^2 + 2\ \Re\left(\K \conj{\U}\right) -2 \abs{\U}^2 = \abs{\K}^2 - \abs{\U}^2 = \abs{ \U^2 -\V} - \abs{\U}^2.
\end{align*}
This expression is clearly independent of $\zeta_{\d}$ and the same is true for $v$.
{The equation for the length interpolation thus simplifies to
$$ \abs{\U^2 - \V} - \abs{\U}^2 +  v = 0.$$
Observing \eqref{def-of-UandV} and \eqref{v} together with \eqref{eq-W-biarc} we see that $\U$, $\V$ and $v$ depend only on the given data and free parameters $\alpha_0$, $\alpha_1$, $\beta_0$ and $\beta_1$.
To emphasize the dependence of the final equation on these parameters, we write it as
\begin{equation}\label{eq-lengthInt1}
e(\alpha_0,\alpha_1,\beta_0,\beta_1) =0 \quad {\rm where} \quad e := \abs{\U^2 - \V} - \abs{\U}^2 +  v.
\end{equation}
Again,  we have one scalar equation for four unknown parameters $\alpha_0$, $\alpha_1$, $\beta_0$ and $\beta_1$.}
To further simplify the  solvability analysis, we additionally assume that
$\beta_0 = \beta_1 = 0$. Since
$$\dr{\r}(0)=\alpha_0^2 \t_0, \; \ddr{\r}(0)= \frac{1}{2} \left(\beta_0 \t_0 + \kappa_0 \alpha_0^4 \n_0\right)
\quad {\rm and} \quad
\dr{\r}(1)=\alpha_1^2 \t_1, \;  \ddr{\r}(1)=\frac{1}{2} \left(\beta_1 \t_1 + \kappa_1 \alpha_1^4 \n_1\right)
$$
this assumption implies that the first and the second order derivative at boundary points are orthogonal; the property that holds for all the points in the case of an arc-length parameterization. Moreover, by introducing a new parameter
$\lambda:=\frac{\abs{\alpha_1}}{\abs{\alpha_0}}$, such that $\lambda^2$ prescribes the ratio between lengths of boundary tangents $\dr{\r}(1)$ and $\dr{\r}(0)$,
we formulate the next theorem.
\begin{thm} \label{thm-2}
Suppose $\r$ is the interpolating planar PH biarc \eqref{def-biarc} of degree $7$, given in Theorem~\ref{thm-main-1} with $\beta_0 = \beta_1 = 0$. 
For any chosen $L>\abs{\Delta \P_0}$ and any positive $\lambda \in \RR$, there exist at least one positive $\alpha^{(1)}_0\in \RR$, such that $e(\pm \alpha^{(1)}_0,  \pm \lambda \alpha^{(1)}_0,0,0) = 0$ and at least one positive $\alpha^{(2)}_0\in \RR$, such that $e(\pm \alpha^{(2)}_0, \mp \lambda \alpha^{(2)}_0,0,0) = 0$.
\end{thm}
\begin{prf}
Since $\t_0$ and $\t_1$ are normalized, we can express them as
$\t_\ell = \cos{\theta_\ell} + \i \sin{\theta_\ell}$, $\ell=0,1$, for some angles $\theta_\ell \in[0,2\pi)$. Under the assumption $\beta_0=\beta_1=0$ and $\alpha_1 = \pm \lambda \alpha_0$, it is easy to see from \eqref{eq-W-biarc} that
$e(\alpha_0, \pm \lambda \alpha_0,0,0)$ includes only even terms of $\alpha_0$, i.e. it is an even function of the unknown $\alpha_0$.
Moreover,
expression {$e_1(\alpha_0, \pm \lambda \alpha_0,0,0)$, $e_1 := -\abs{\U}^2 + v$,}
is an even polynomial of degree $6$ in a variable $\alpha_0$, and
{$e_2(\alpha_0, \pm \lambda \alpha_0,0,0)$, $e_2 := |\U^2 - \V|^2$,}
 is an even
polynomial of degree $12$ in a variable $\alpha_0$. We further compute 
that the free and the leading coefficient of
{$e_1(\alpha_0, \pm \lambda \alpha_0,0,0)$}
are equal to
$$
c_{0} = -\frac{140}{13} L,\quad c_{1}=\frac{1}{29952}\left(131 \kappa_0^2 \pm 61 \Theta \lambda^3 \kappa_0 \kappa_1 +
131 \lambda^6 \kappa_1^2\right)$$
where
$$
\Theta = \sqrt{\cos \left(\theta _0\right)+1} \sqrt{\cos \left(\theta _1\right)+1} \left(\tan \left(\frac{\theta _0}{2}\right) \tan
   \left(\frac{\theta _1}{2}\right)+1\right),
$$
and that $
c_{2} = (\frac{140}{13})^2 \abs{\Delta \P_0}^2$ is the free coefficient of
{$e_2(\alpha_0, \pm \lambda \alpha_0,0,0)$}
Thus, as $e= \sqrt{e_2}+e_1$, we get $$e(0,0,0,0) = \frac{140}{13}\left(\abs{\Delta \P_0} - L\right) < 0.$$
If we prove that $c_{1}$ is always positive, then the limit $\displaystyle{\lim_{\alpha_0 \to \infty} e(\alpha_0, \pm \lambda \alpha_0,0,0)}$ is also positive. Thus there exist at least one positive zero $\alpha^{(1)}_0$ of
the function $e(\alpha_0,\lambda \alpha_0,0,0)$, and at least one positive zero $\alpha^{(2)}_0$ of the function
$e(\alpha_0,-\lambda \alpha_0,0,0)$. Since these two functions are even, the statement of the theorem follows.

It remains to show that $c_1 > 0$ independently of the chosen $\lambda$, $\kappa_\ell$ and $\theta_\ell$, $\ell=0,1$.
From trigonometric identities it follows that
$\abs{\Theta} = 
2\abs{\cos{\left(\frac{1}{2}(\theta_0 - \theta_1)\right)}}\leq 2,
$
and this implies
\begin{align*}
29952 c_1 =&  131 \kappa_0^2 \pm 61 \Theta \lambda^3 \kappa_0 \kappa_1 + 131 \lambda^6 \kappa_1^2 \geq
131 \kappa_0^2 - 122 \lambda^3 \abs{\kappa_0 \kappa_1} + 131 \lambda^6 \kappa_1^2 \\
&> 131 \kappa_0^2 - 262 \lambda^3 \abs{\kappa_0 \kappa_1} + 131 \lambda^6 \kappa_1^2 = \left(\sqrt{131} \abs{\kappa_0} - \sqrt{131} \lambda^3 \abs{\kappa_1}\right)^2 \geq 0,
\end{align*}
which concludes the proof.
\qed
\end{prf}

From Theorem~\ref{thm-2} it follows that there exist at least four solutions of the equation \eqref{eq-lengthInt1} for $\beta_0=\beta_1=0$ and a fixed parameter $\lambda$.  Since this equation is independent of the sign $\zeta_{\d}$, this would imply at least eight biarc solutions, obtained for two different choices of $\zeta_{\d} \in \{-1,1\}$. However, simultaneously  changing the signs of $\alpha_0$, $\alpha_1$ and $\zeta_{\d}$ does not change the control points of the biarc. Therefore, we fix the sign $\zeta_{\d}=1$, and in general we have four solution interpolants.


\section{Examples}   \label{sec:Examples}

From what we have discussed in the previous section, the $G^2$ interpolating PH biarc fulfilling also the prescribed length constraint can be computed by solving the nonlinear equation $e(\alpha_0, \lambda \alpha_0,\beta_0,\beta_1) =0$, given in \eqref{eq-lengthInt1}, for the unknown $\alpha_0$. This equation depends on three free parameters $\beta_0$, $\beta_1$ and $\lambda$. {By fixing $\beta_0, \beta_1$ and $\lambda$ the problem reduces to the solution of a non-linear (algebraic) equation in one unknown, which can be computed easily by a few Newton-Raphson iterations.} Examples in this section show that the choice $\beta_0 = \beta_1 = 0$ and $\lambda=1$  leads to interpolants with nice shape properties. More precisely, the fairness of the resulting curves is measured by considering the {\it bending energy}
 \be\label{bend_en}
 E=\int_0^1 \kappa^2(s) \ ds
 \ee
associated to the biarc curve (see for instance \cite{FPS'21}).
{
For data far from being symmetric, we have also investigated the case where $\lambda$ is obtained by minimizing the bending energy under the length constraint given by \eqref{eq-lengthInt1} and the assumption $\beta_0 = \beta_1 = 0$.
The optimization procedure has been implemented in {\tt Matlab} by using the solver routines in its Optimization toolbox. In more detail we have used {\tt fmincon} which uses the Interior Point algorithm. We note that, despite of the results of Theorem \ref{thm-2}, a reasonable solution could be obtained also by minimizing \eqref{bend_en} with respect to all the free parameters $\beta_0$, $\beta_1$ and $\lambda$, as we have shown in the Example~\ref{Sec5-example1}.
}

Without lack of generality, we assume in the first four examples that the $G^1$ data - points and tangent directions - are chosen as
in \eqref{example-int-data}. In the first example we 
give a comparison between the derived PH biarc interpolants and the single PH curve interpolants, considered in Section~\ref{sec:probl}.

\begin{exm}  \label{Sec5-example1}
Let the data be chosen as in Example~\ref{exm-OnePHCurve}, i.e. with
$\theta_0=-\frac{\pi}{4}$, $\theta_1=-\frac{\pi}{8}$, $\kappa_0 = 1$ and $\kappa_1 = -1$.
For both prescribed lengths, $L=1.1$ and $L=1.05$, and the free parameters chosen as $\beta_0=\beta_1=0$, $\lambda=1$, there exist four different biarc solutions, described in Theorem~\ref{thm-main-1}, with $\zeta_\d=1$ and
$\alpha_0$, $\alpha_1$ given in Table~\ref{Tab:example-1}. The last column of this table shows values of the corresponding bending energy \eqref{bend_en} for each solution, and the biarc having this minimal value is shown in Figure~\ref{Fig:ex_5_1a} -- left for $L=1.1$, right for
$L=1.05$. Comparing this `best' PH biarc (for $L=1.1$) with the two PH interpolants from Example~\ref{exm-OnePHCurve} - shown in Figure~\ref{fig:example-int-data}, we see that the biarc construction provides visually better results, comparing the shape of the control polygon as well as the porcupine curvature plot. It also performs better regarding the bending energy that equals $E=184.113$, for the left and $E=14.4481$ for the right PH curve from Figure~\ref{fig:example-int-data}.
\begin{table}
\begin{center}
\renewcommand\arraystretch{1.1}
\begin{tabular}{|c|c|c|c|c|c||r|r|r|}
\hline
 & $\theta_0$ & $\theta_1$ &  $\kappa_0$ & $\kappa_1$ & $L$ & $\alpha_0 \quad $ & $\alpha_1 \quad $ & Bending energy \\
\hline
I & \multirow{4}*{$-\frac{\pi}{4}$} & \multirow{4}*{$ -\frac{\pi}{8}$} &
\multirow{4}*{$1$} & \multirow{4}*{$-1$} & \multirow{4}*{$1.1$} &
$1.15932$  & $1.15932$ & $6.01964\cdot 10^0$\\
\cline{7-9}
II & & & & &  & $-1.15932$  & $-1.15932$ & $1.64506\cdot 10^6$\\
\cline{7-9}
III & & & & &  & $-0.96713$  & $0.96713$ & $1.03930 \cdot 10^4$\\
\cline{7-9}
IV & & & & &  & $0.96713$  & $-0.96713$ & $3.44494 \cdot 10^4$\\
 \hline
 \hline
 I  & \multirow{4}*{$-\frac{\pi}{4}$} & \multirow{4}*{$ -\frac{\pi}{8}$} &
\multirow{4}*{$1$} & \multirow{4}*{$-1$} & \multirow{4}*{$1.05$} &
$0.85919$  & $0.85919$ & $4.85785\cdot 10^0$\\
\cline{7-9}
II & & & & &  & $-0.85919$  & $-0.85919$ & $4.60342\cdot 10^6$\\
\cline{7-9}
III & & & & &  & $-0.72422$  & $0.72422$ & $3.81363 \cdot 10^4$ \\
\cline{7-9}
IV & & & & &  & $0.72422$  & $-0.72422$ & $2.41983 \cdot 10^5$\\
 \hline
 \hline
 I  & \multirow{4}*{$-\frac{\pi}{4}$} & \multirow{4}*{$ \frac{\pi}{8}$} &
\multirow{4}*{$1$} & \multirow{4}*{$1$} & \multirow{4}*{$1.1$} &
$1.31430$  & $1.31430$ & $5.15473\cdot 10^0$\\
\cline{7-9}
II & & & & &  & $-1.31430$  & $-1.31430$ & $1.70151\cdot 10^7$\\
\cline{7-9}
III & & & & &  & $-1.24343$  & $1.24343$ & $6.04317 \cdot 10^4$ \\
\cline{7-9}
IV & & & & &  & $1.24343$  & $-1.24343$ & $1.40778 \cdot 10^5$\\
 \hline
 \hline
 \end{tabular}
\end{center}
\caption{Data values, resulting $\alpha_0$, $\alpha_1$, and the corresponding bending energy of the four PH biarcs from Example~\ref{Sec5-example1}.  }
\label{Tab:example-1}
\end{table}
\begin{figure}[htb]
\centering
\includegraphics[width=0.48\textwidth]{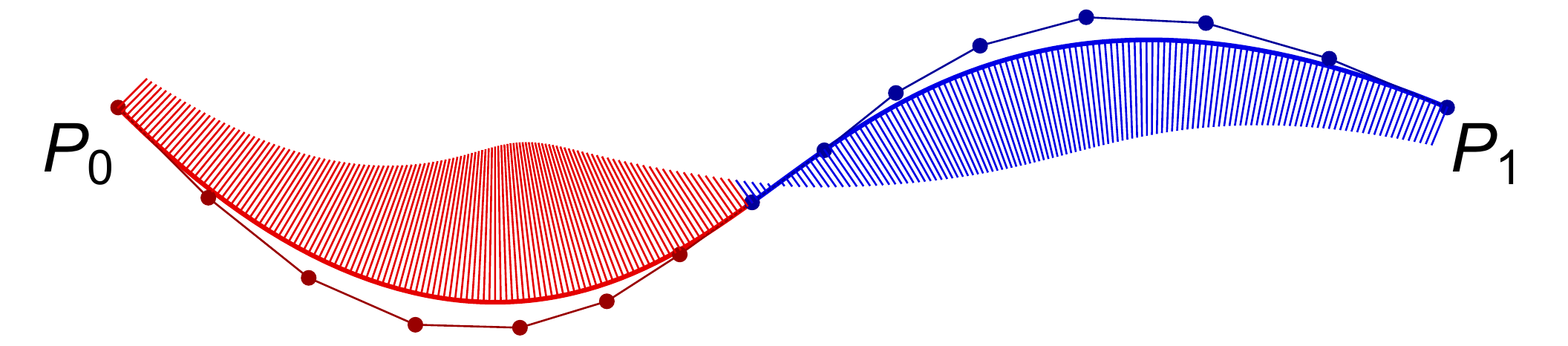}
\hskip.5cm
\includegraphics[width=0.48\textwidth]{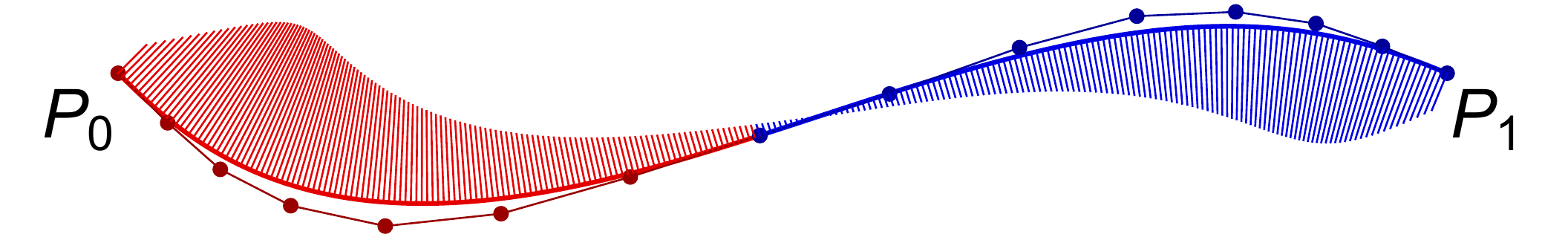}
\caption{The PH biarc interpolant with assigned length $L$ of minimal bending energy, computed by choosing $\beta_0=\beta_1=0$ and $\lambda=1$, together with the control polygon and the porcupine curvature plot (with proportional factor equal to 0.03). Left: $L=1.1$; right:$L=1.05$.
}
\label{Fig:ex_5_1a}
\end{figure}
\begin{figure}[htb]
\centering
\begin{minipage}{0.48\textwidth}
\includegraphics[width=1\textwidth]{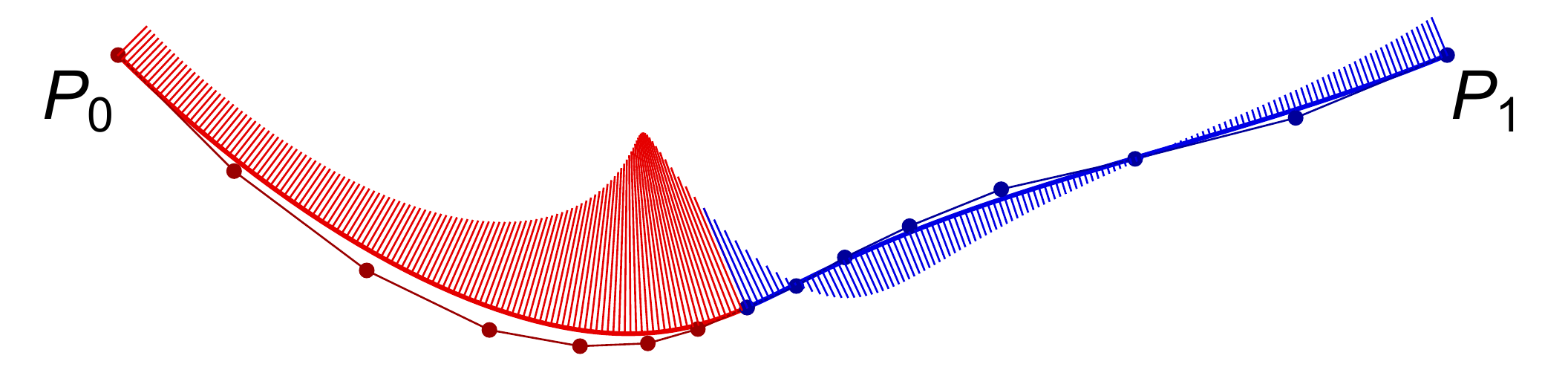}
\end{minipage}
\hskip.5cm
\begin{minipage}{0.48\textwidth}
\includegraphics[width=1\textwidth]{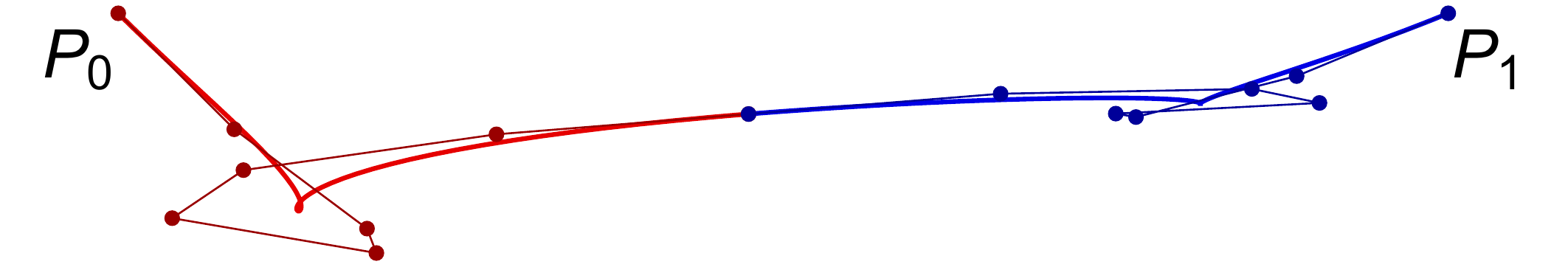}
\end{minipage}\\
\begin{minipage}{0.48\textwidth}
\includegraphics[width=1\textwidth]{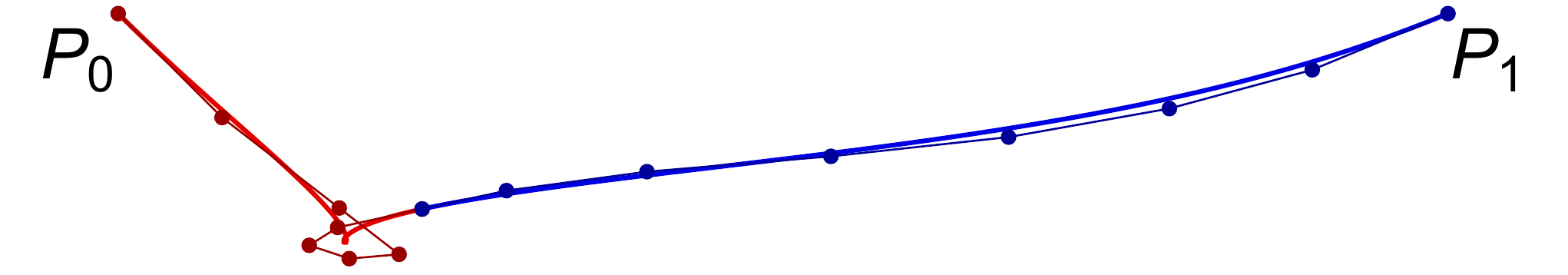}
\end{minipage}
\hskip.5cm
\begin{minipage}{0.48\textwidth}
\includegraphics[width=1\textwidth]{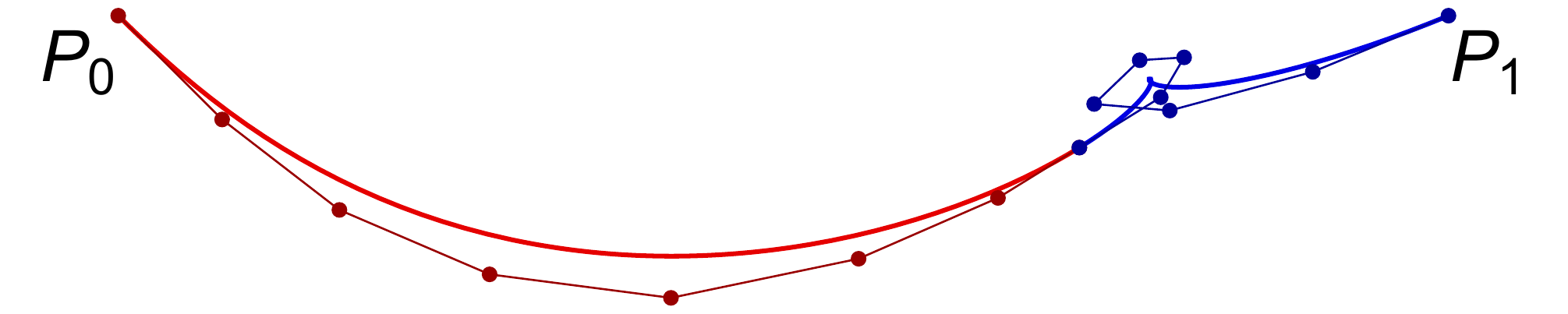}
\end{minipage}
\caption{Four PH biarc interpolants of degree $7$ for the data \eqref{example-int-data} with
$\theta_0 = -\frac{\pi}{4}$, $\theta_1 = \frac{\pi}{8}$, $\kappa_0 = \kappa_1 = 1$, $L=1.1$, and
the free parameters chosen as $\beta_0=\beta_1 =0$ and $\lambda=1$. The graph of the biarc with the minimal bending energy is equipped with the porcupine curvature plot  (with proportional factor equal to 0.03).
}
\label{Fig:ex_5_1b}
\end{figure}

Figure~\ref{Fig:ex_5_1b} shows all four biarc solutions corresponding to the third set (rows $9$--$12$) of Table~\ref{Tab:example-1}. These curves can be compared to the ones shown in Figure~\ref{fig:example-int-data-a}, where the single PH curve interpolants for the same data values are considered. We note that in that case the bending energies equal to $1.02189\cdot 10^6$, $183.06$, $2.06226\cdot 10^6$, $4.46632\cdot 10^5$, respectively. Again, the biarcs perform much better.
We could further improve the solution by optimizing $E$ with respect to the parameter $\lambda$.
Figure~\ref{Fig:ex_5_1c} (top left) shows values of $E$ (in $\log_{10}$ scale) in dependence of $\lambda = \frac{j}{10}$, $j=1,2,\dots,100$ (with $\beta_0=\beta_1 = 0$), where among different solutions
(for fixed values of free parameters) we choose the one with the minimal value of $E$.
The minimum is attained at $\lambda=0.5$ and equals $E=1.559197$. 
On Figure~\ref{Fig:ex_5_1c} (top right), the PH biarcs are shown for $\lambda \in \{0.1, 0.5, 1, 1.5,2,3,4,5,10\}$. {Of course, this discrete optimization can be replaced by the constrained optimization solver to compute the value of the optimal parameter $\lambda$ even more precisely. Namely, we get that the minimal value of $E$ (for fixed $\beta_0=\beta_1=0$) equals $1.553895$ and it is attained at
$\lambda=0.521524$. However, the difference between both minimal values is negligibly small.}
\begin{figure}[htb]
\centering
\begin{minipage}{0.45\textwidth}
\includegraphics[width=1\textwidth]{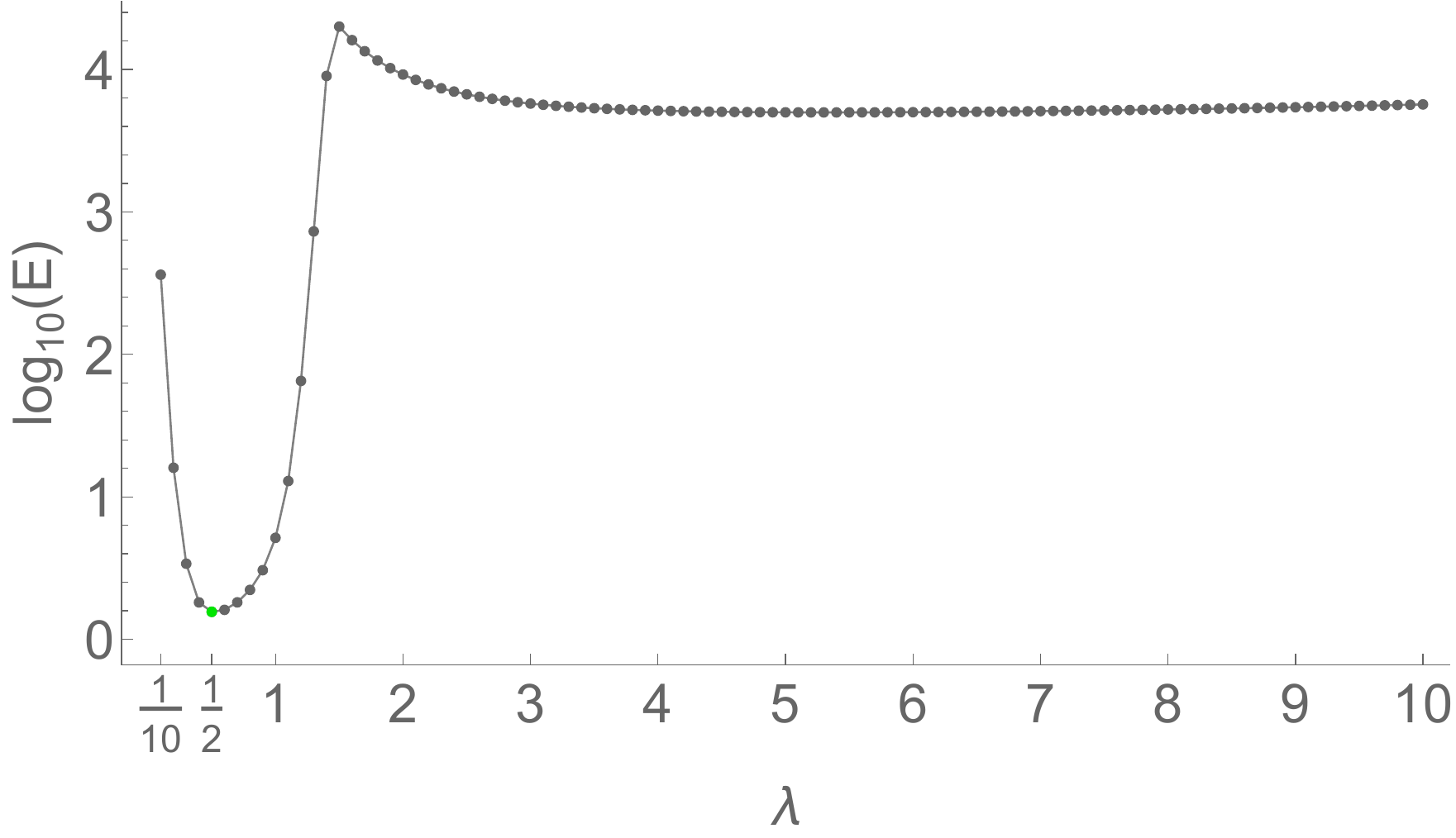}
\end{minipage}
\hskip.5cm
\begin{minipage}{0.45\textwidth}
\includegraphics[width=1\textwidth]{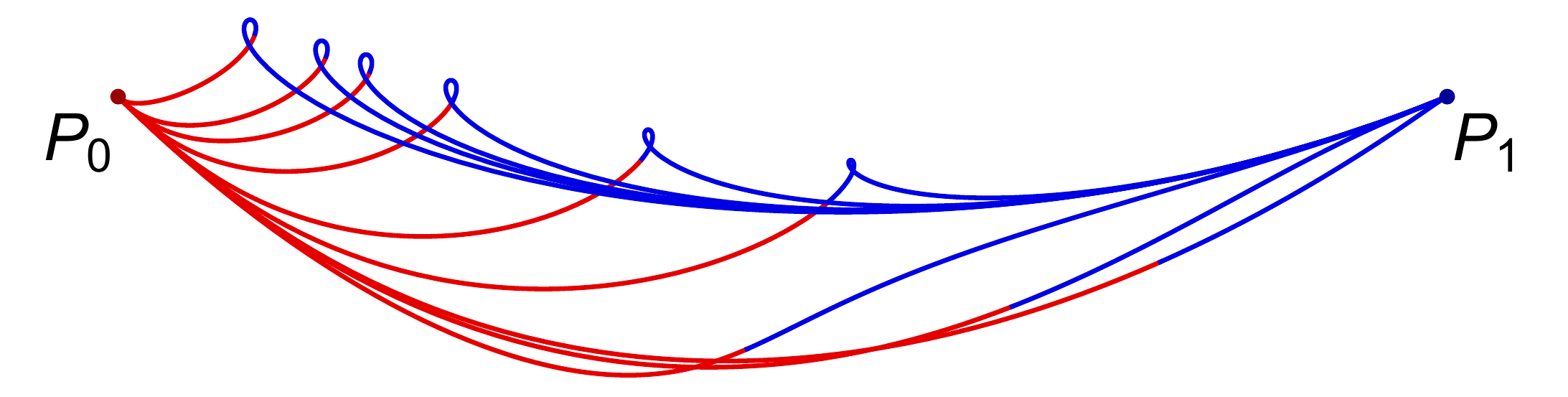}
\end{minipage}\\
\begin{minipage}{0.45\textwidth}
\includegraphics[width=1\textwidth]{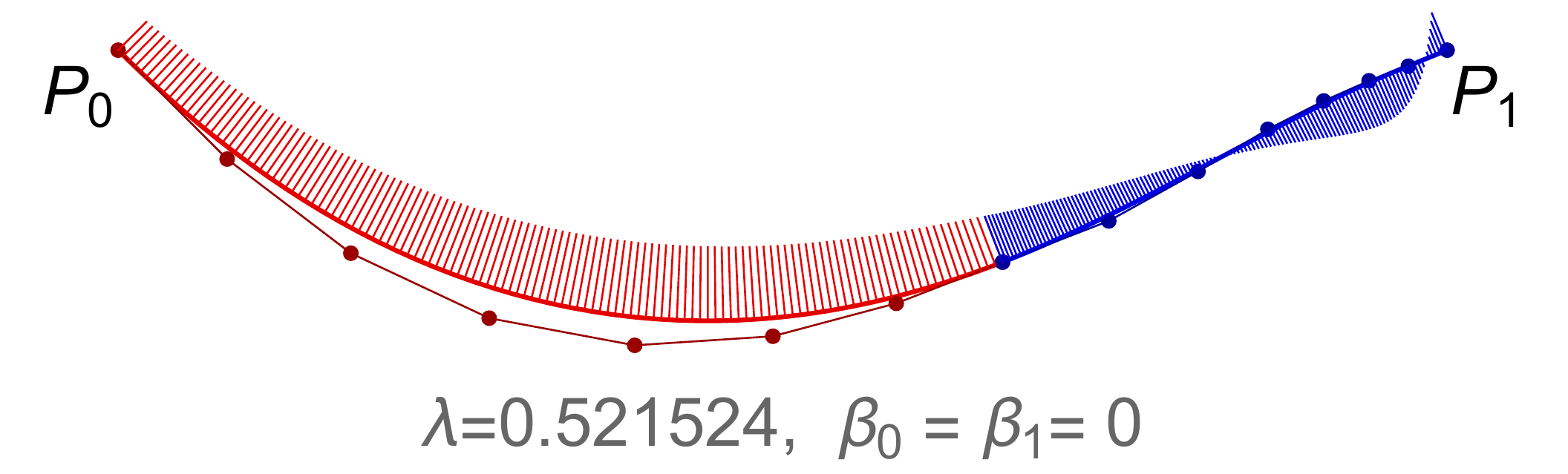}
\end{minipage}
\hskip.5cm
\begin{minipage}{0.45\textwidth}
\includegraphics[width=1\textwidth]{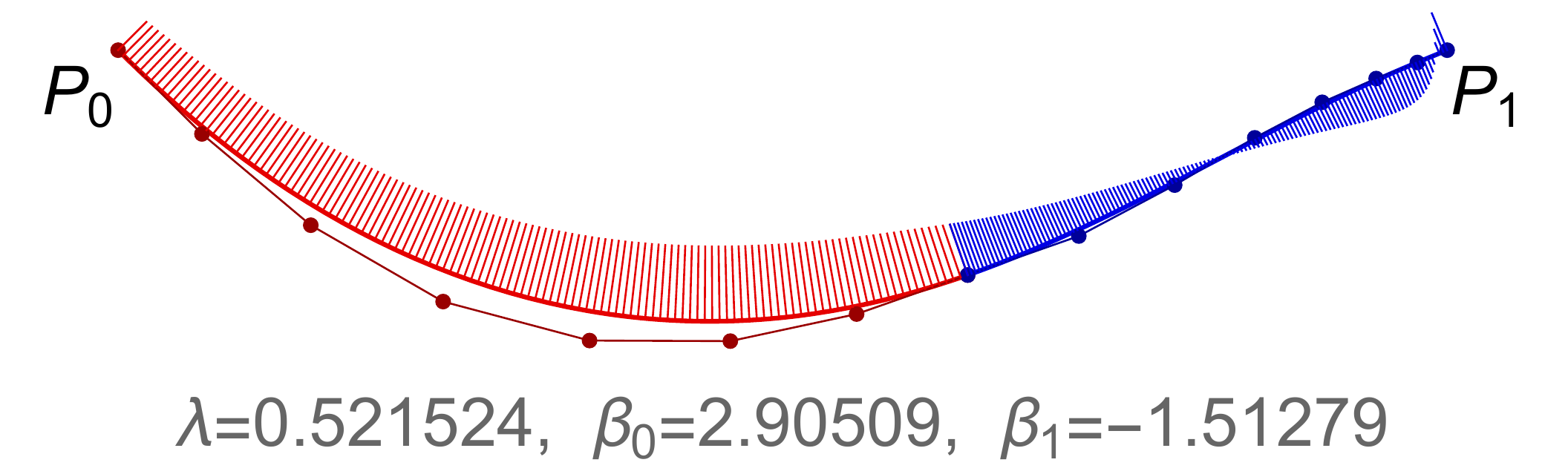}
\end{minipage}
\caption{
Top: The bending energy $E$ and PH biarcs in dependence of the free parameter $\lambda$ for the data \eqref{example-int-data} with $\theta_0 = -\frac{\pi}{4}$, $\theta_1 = \frac{\pi}{8}$, $\kappa_0 = \kappa_1 = 1$, $L=1.1$, and
$\beta_0=\beta_1 =0$. Bottom: Two optimal PH biarcs.
}
\label{Fig:ex_5_1c}
\end{figure}
Using the constrained optimization solver we can minimize also with respect to $\beta_0$ and $\beta_1$, although theoretically the existence of the solution if these two parameters are non-zero is not guaranteed.
{Choosing $\lambda=0.521524$ we have that the minimal value of $E$ is achieved for $\beta_0 = 2.90509$ and
$\beta_1 = -1.51279$, with $E=1.543123$.} However, comparing this optimal PH biarc with the one computed by
{$ \lambda=0.521524$,} $\beta_0=\beta_1=0$, we see that the difference in the value of $E$ as well as on their graphs (see Figure~\ref{Fig:ex_5_1c}) is negligible.
\end{exm}

The data in the following three examples are taken from \cite{Farouki'16}.
In all the cases we fix $\beta_0=\beta_1=0$.
\begin{exm}  \label{Sec5-example2}
Let us consider an example of parallel end tangents, obtained by taking $\theta_0=\theta_1=\pi/4$. Let us then assign
$\kappa_0=-0.5$, $\kappa_1=0.5$ and $L=1.5$. In this case it is reasonable to set $\lambda=1$.
All four solutions 
are reported in Table~\ref{Tab:example-2} and shown in Figure~\ref{Fig:ex_5_2}. The one with the minimal bending energy is equipped with the porcupine curvature plot.
\begin{table}[htb]
\begin{center}
\renewcommand\arraystretch{1.1}
\begin{tabular}{|c|c|c|c|c|c||r|r|c|}
\hline
 & $\theta_0$ & $\theta_1$ &  $\kappa_0$ & $\kappa_1$ & $L$ & $\alpha_0 \quad$ & $\alpha_1 \quad$ & Bending energy \\
\hline
I  & \multirow{4}*{$\frac{\pi}{4}$} & \multirow{4}*{$ \frac{\pi}{4}$} &
\multirow{4}*{$-0.5$} & \multirow{4}*{$0.5$} & \multirow{4}*{$1.5$} &
$1.60884$  & $1.60884$ & $26.2939$\\
\cline{7-9}
II & & & & &  & $-1.60884$  & $-1.60884$ & $2928.06$\\
\cline{7-9}
III & & & & &  & $-1.31667$  & $1.31667$ & $239.358$\\
\cline{7-9}
IV & & & & &  & $1.31667$  & $-1.31667$ & $239.358$\\
 \hline
 \end{tabular}
\end{center}
\caption{Data values, solutions $\alpha_0$, $\alpha_1$, and the bending energy of the PH biarcs from Example~\ref{Sec5-example2}.  }
\label{Tab:example-2}
\end{table}
\begin{figure}[htb]
\centering
\begin{minipage}{0.48\textwidth}
\includegraphics[width=0.85\textwidth]{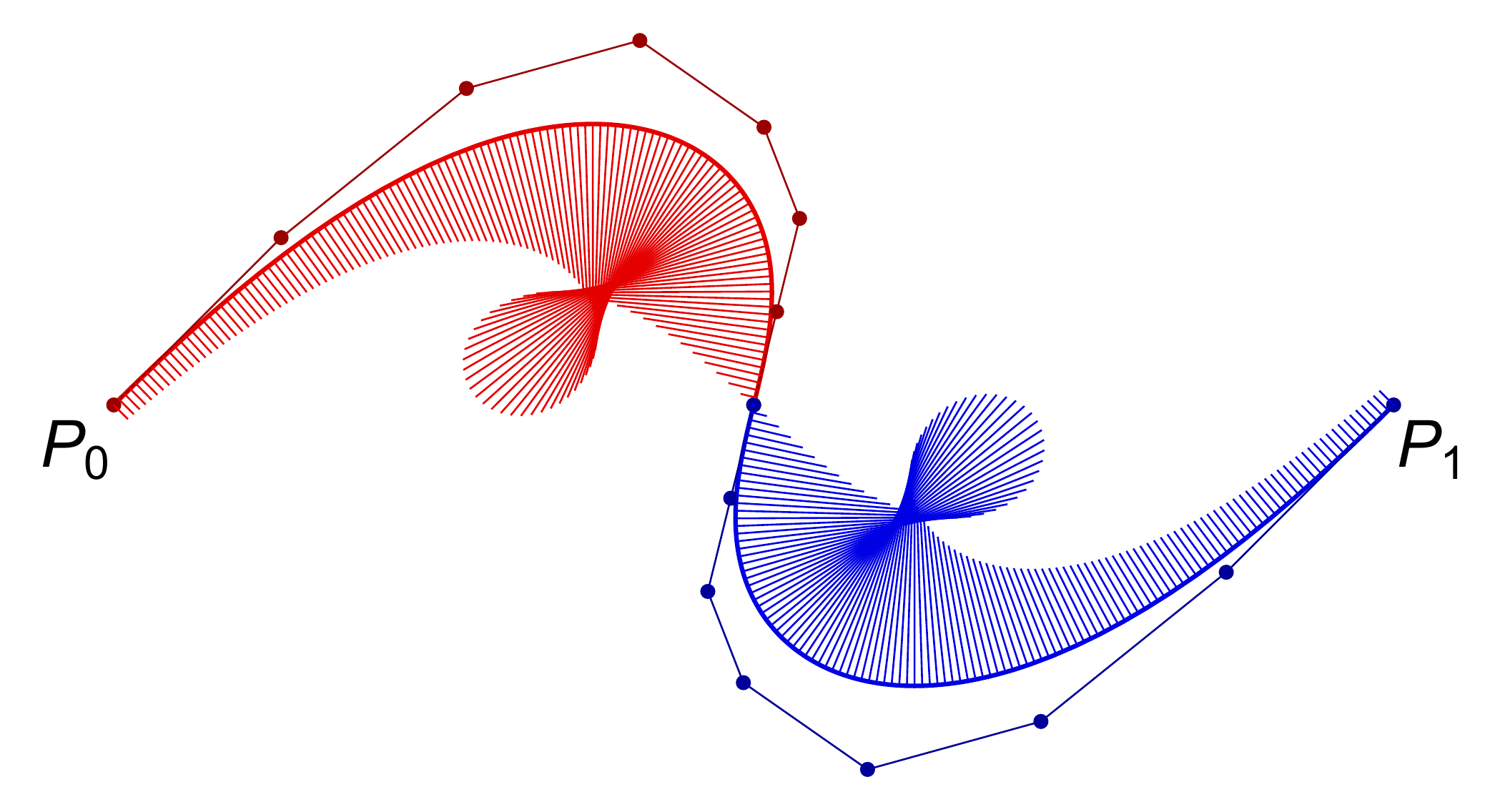}
\end{minipage}
\hskip.5cm
\begin{minipage}{0.48\textwidth}
\includegraphics[width=0.85\textwidth]{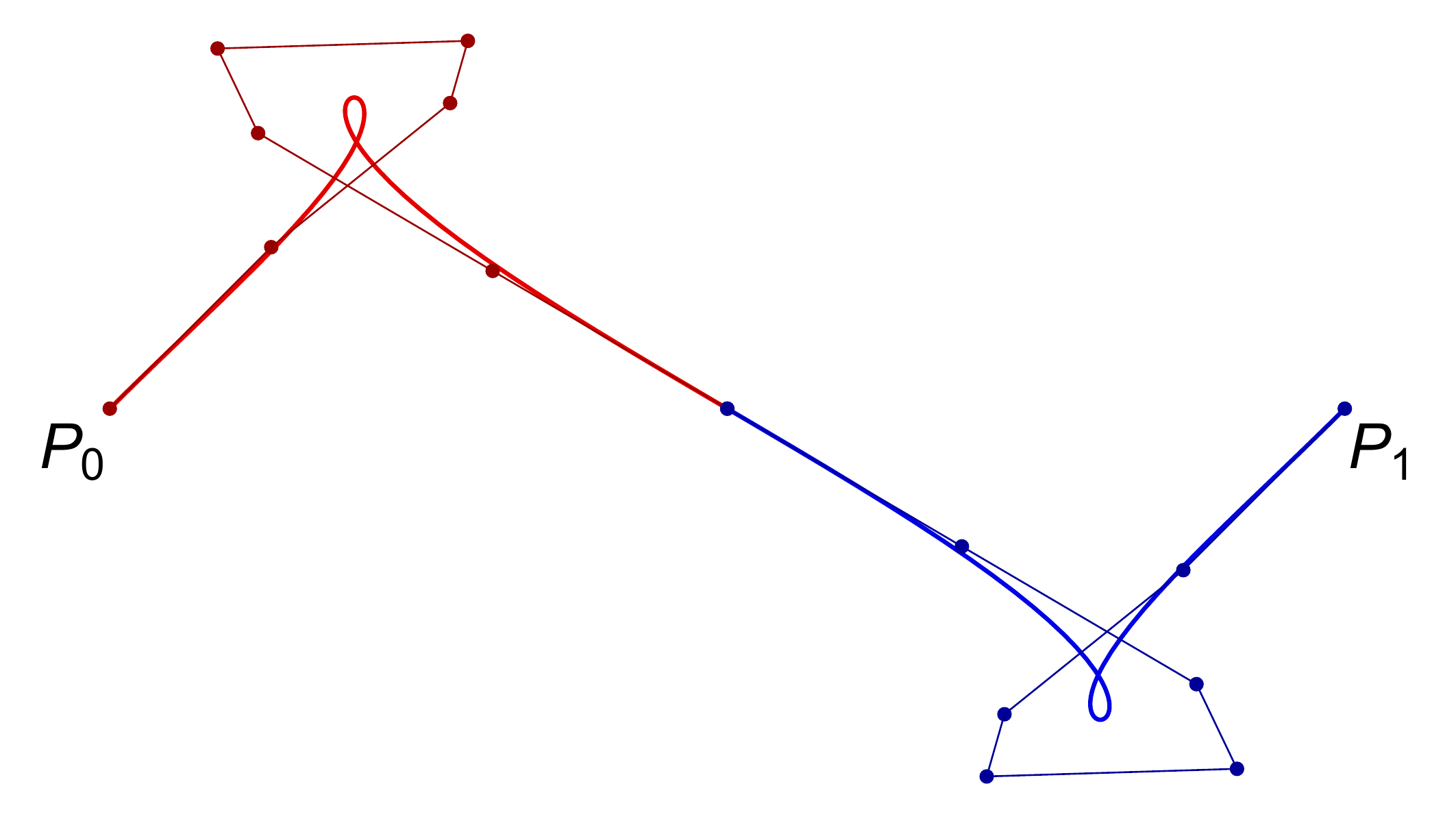}
\end{minipage}\\
\begin{minipage}{0.48\textwidth}
\includegraphics[width=0.85\textwidth]{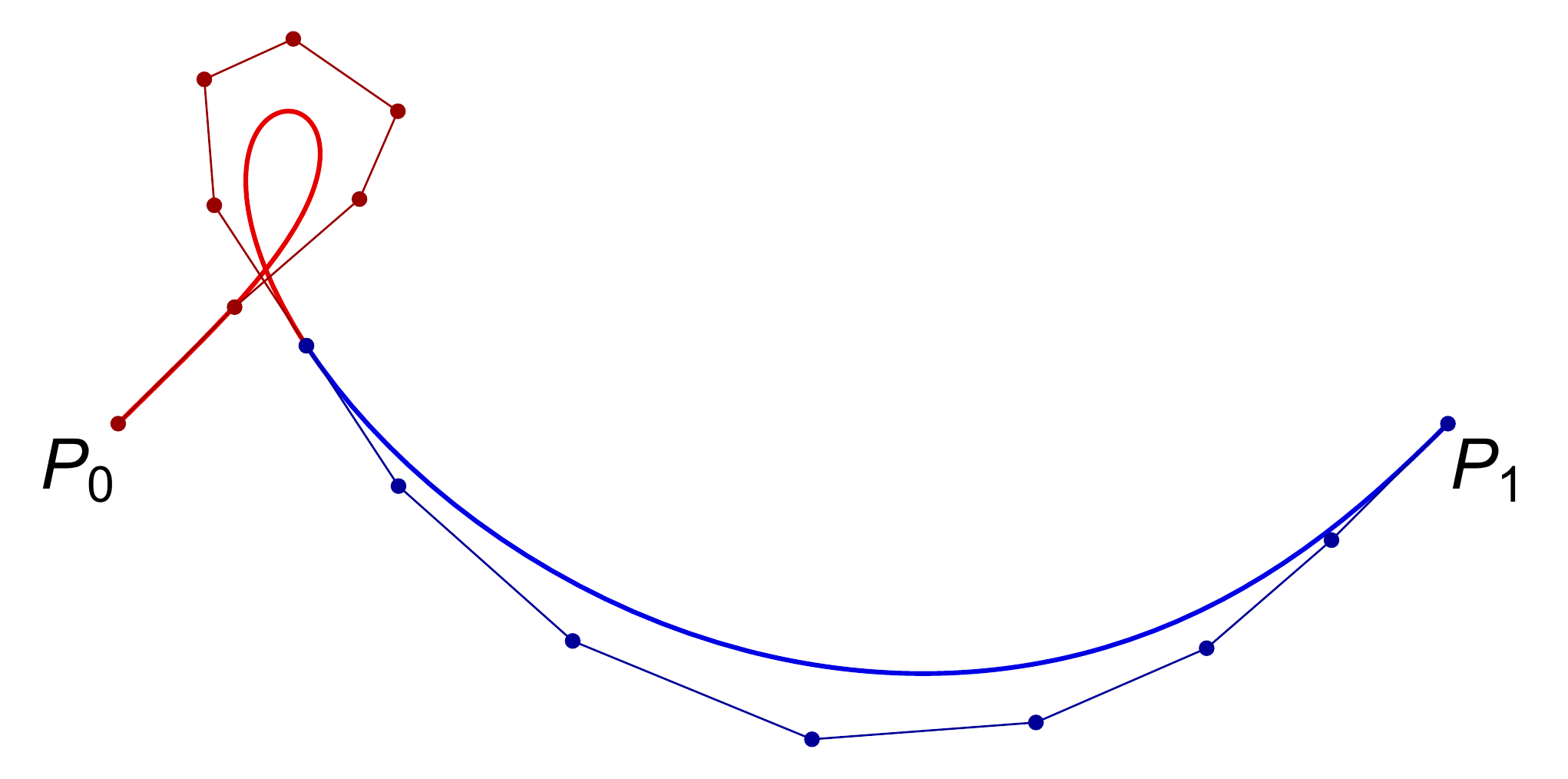}
\end{minipage}
\hskip.5cm
\begin{minipage}{0.48\textwidth}
\includegraphics[width=0.85\textwidth]{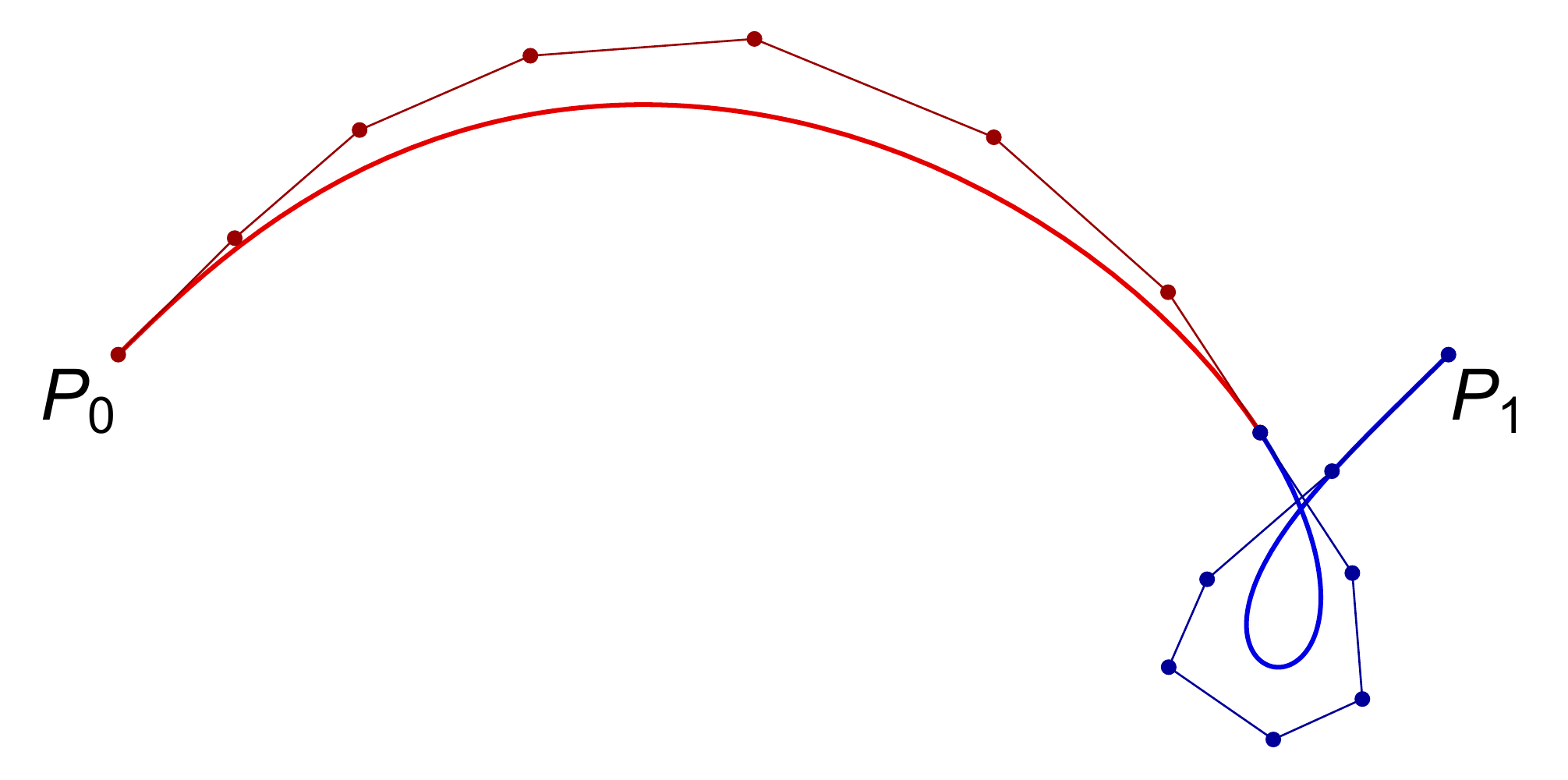}
\end{minipage}
\caption{Four PH biarc interpolants of degree $7$ for the data \eqref{example-int-data} with
$\theta_0 = \theta_1= \frac{\pi}{4}$,  $\kappa_0 = -\kappa_1 = -0.5$, $L=1.5$, and
the free parameters chosen as $\beta_0=\beta_1 =0$ and $\lambda=1$. The graph of the biarc with the minimal bending energy is equipped with the porcupine curvature plot.
}
\label{Fig:ex_5_2}
\end{figure}
\end{exm}

\begin{exm}  \label{Sec5-example3}
In the third example we consider symmetric data, where $\theta_0=-\theta_1=\pi/3$, $L=1.35$ and $\kappa_0=\kappa_1=-0.5$. Also in this case, it is reasonable to set $\lambda=1$.
 The resulting curve with the minimal bending energy is shown in Figure~\ref{Fig:ex_3}, along with the porcupine curvature plot (with proportional factor equal to 0.03). The graph of the curvature is also shown. We observe that the symmetry of the data is preserved and the length of both biarc segments is equal to $L/2$. In this case we have $\alpha_0=\alpha_1=1.27991$, while the bending energy is $3.51446$.
\begin{figure}[htb]
\centering
\begin{minipage}{0.53\textwidth}
\includegraphics[width=1\textwidth]{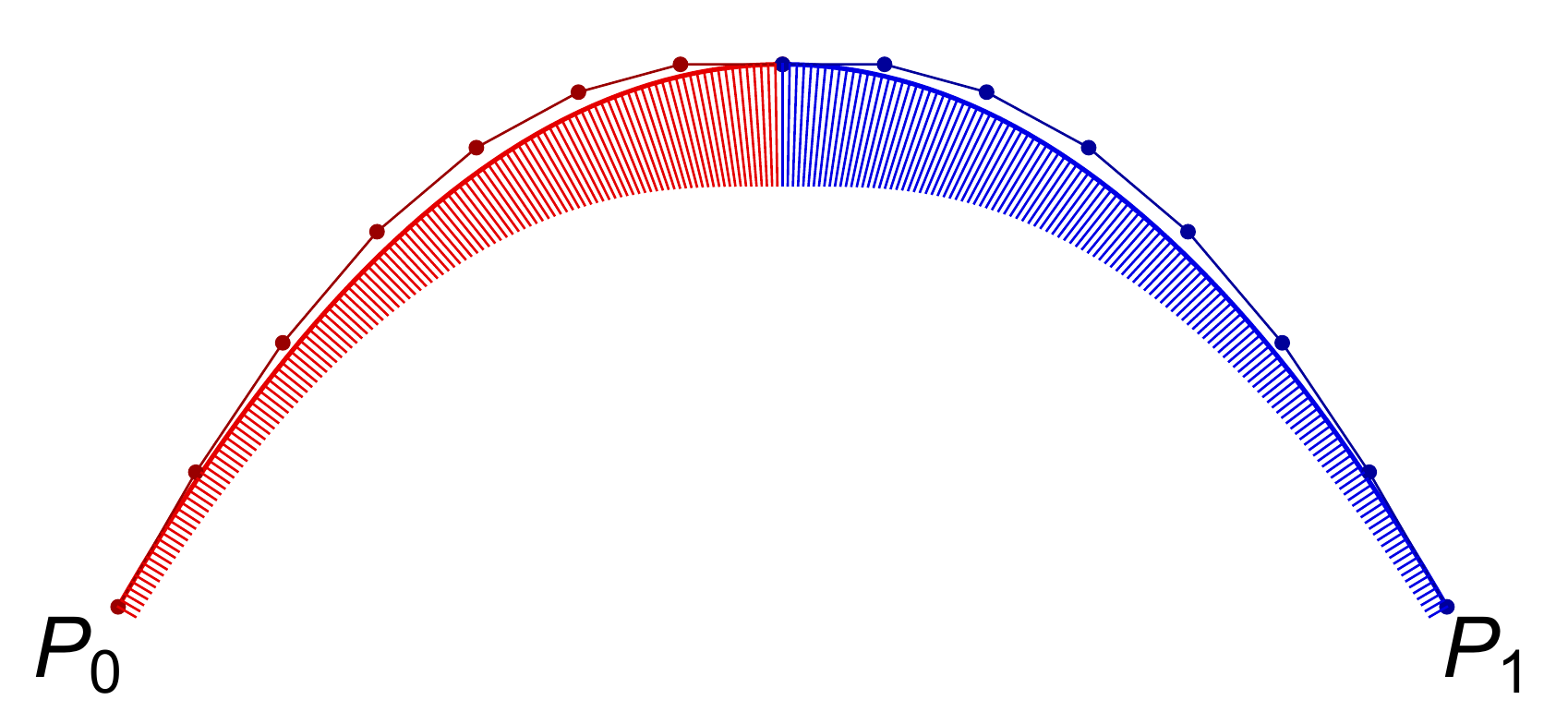}
\end{minipage}
\hskip.5cm
\begin{minipage}{0.42\textwidth}
\includegraphics[width=1\textwidth]{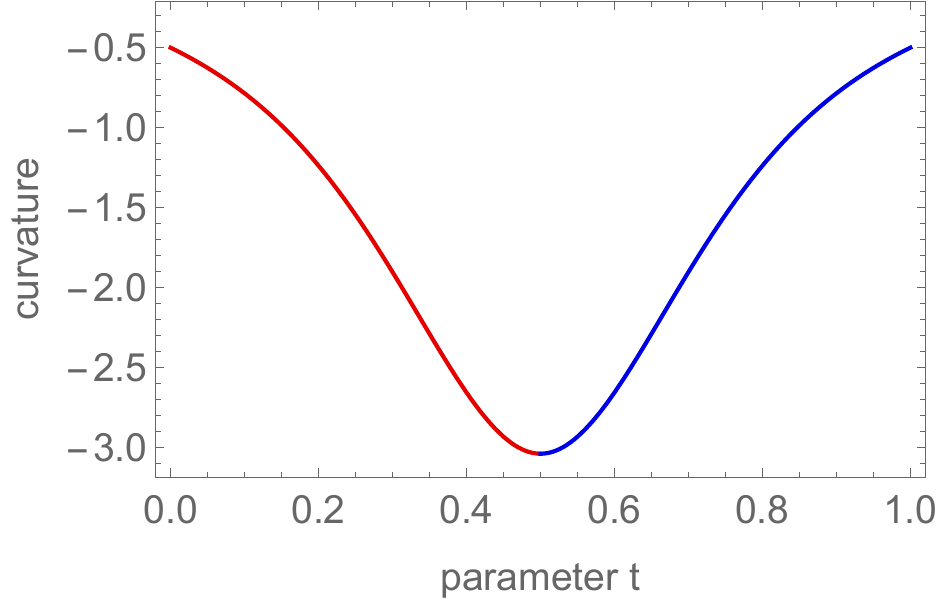}
\end{minipage}
\caption{
Example~\ref{Sec5-example3}. The resulting curve along with the porcupine curvature plot (left), and the graph of the curvature (right).
}
\label{Fig:ex_3}
\end{figure}

\end{exm}
\begin{exm}  \label{Sec5-example4}
We now consider two sets of data: a convex one with $\theta_0=\pi/4$ and $\theta_1=-\pi/3$, and a non-convex data set obtained by taking $\theta_0=-\pi/6,$ and $\theta_1=-\pi/3$. In Figure~\ref{Fig:ex_4} we show the resulting curves for different lengths, ranging from $L=1.1$ to $L=1.6$. For simplicity the assigned curvatures $\kappa_0$ and $\kappa_1$ are always set to zero, and the free parameter $\lambda$ is set to one.
Among all the solutions, we pick the one with the minimal bending energy, obtained (in all the cases)
from positive solutions for $\alpha_0$ and $\alpha_1$.
\begin{figure}[htb]
\centering
\begin{minipage}{0.45\textwidth}
\includegraphics[width=1\textwidth]{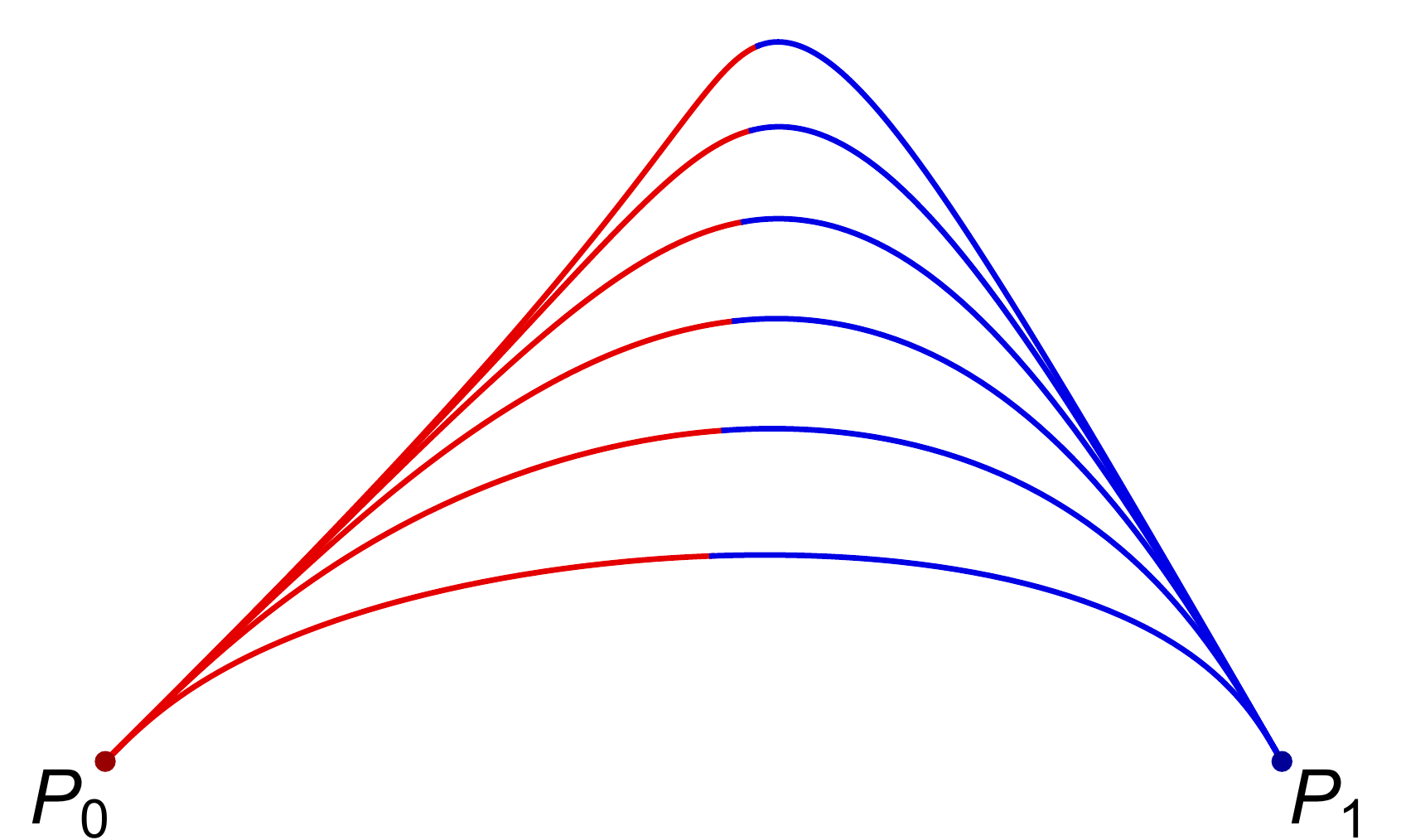}
\end{minipage}
\hskip.5cm
\begin{minipage}{0.45\textwidth}
\includegraphics[width=1\textwidth]{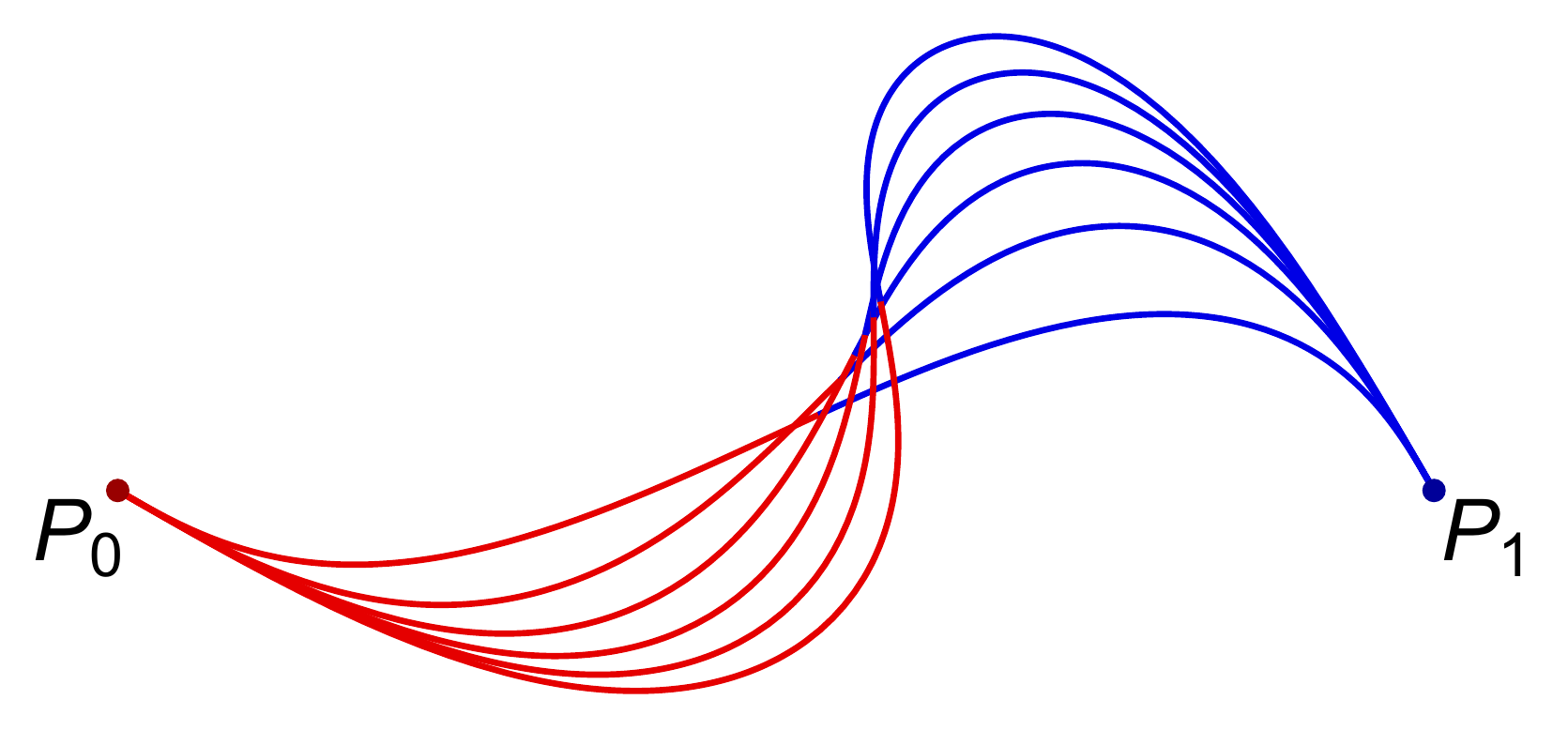}
\end{minipage}
\caption{
Example~\ref{Sec5-example4}. The resulting curves with minimal bending energy for convex and non-convex data with lengths varying from $L=1.1$ to $L=1.6$.
}
\label{Fig:ex_4}
\end{figure}
\end{exm}

Next let us demonstrate the performance of the presented interpolation scheme for curve approximation and its generalization to the spline construction.
\begin{exm}  \label{Sec5-example5}
Suppose that we sample the data from a well known curve - the logarithmic spiral - given by the parameterization (in the complex plane)
\begin{equation}\label{eq:log-spiral}
\f(s) = -e^{\omega s} \cos{(s)}  +   \i \,e^{\omega s} \sin{(s)}, \quad s\in [a,b],
\end{equation}
for some real number $\omega$.
It is straightforward to compute the tangents, the curvature $\kappa_\f$ and the length of this curve from $\f(s_i)$ to $\f(s_f)$:
\begin{align*}
& \f'(s) = e^{s \omega } (\sin (s) - \omega  \cos (s)) + \i \, e^{s \omega } (\cos (s)+ \omega  \sin (s)),\\
& \kappa_\f(s) = -\frac{e^{-s \omega }}{\sqrt{1+\omega ^2}},\quad L_\f(s_i,s_f) =
\frac{\sqrt{1+\omega ^2}}{\omega} \left(e^{s_f \,\omega }-e^{s_i \, \omega }\right).
\end{align*}
First, we choose $a=0$ and $b=h$, for decreasing values $h = 2^{-k}$, $k=0,1,\dots, 8$,  and observe the error
between $\f$ and the PH biarc $\r$, with interpolation data \eqref{eq:intCond} chosen as
\begin{equation} \label{data-from-f}
\P_0 = \f(0),  \; \t_0=\f'(0), \; \kappa_0 = \kappa_\f(0), \quad
\P_1 = \f(h),  \; \t_1=\f'(h), \; \kappa_1 = \kappa_\f(h), \quad L=L_\f(0,h).
\end{equation}
The free parameters are again set to $\beta_0=\beta_1=0$, $\lambda=1$, and among all the solutions we take the one with the minimal bending energy.
To measure the error we choose
\begin{equation} \label{eq:error}
E_{{\rm err}}(\r, \f) = \max_{t\in [0,1]}
\nrm{\r(t) - (\f \circ\varphi)(t)},
\end{equation}
which gives the upper bound for the parametric distance. Here $\varphi:[0,1]\to [a,b]$ is the (bijective) reparameterization function, computed in such a way that $\varphi'>0$ and  $\r^{(\ell)}(j)=\frac{d^\ell \left(\f \circ \varphi\right)}{d t^\ell}(j)$ for $\ell=0,1,2$ and $j=0,1$.
For the choice $\omega=0.2$ the computed errors are given in Table~\ref{Tab:example-5-5}, together with the {\it decay exponent}, and graphically represented in Figure~\ref{Fig:ex_5a}. These results numerically confirm that the PH biarc approximates the given curve with the approximation order $5$. For the sake of comparison, the right part of Table~\ref{Tab:example-5-5} shows the errors of the single PH curve interpolant from Section~\ref{sec:probl}. In this case the order of approximation is $6$, which is the expected optimal order when interpolating $G^2$ data. As usual, the order of approximation decreases by one when the biarcs are used, but the interpolation scheme becomes much simpler.
\begin{table}[htb]
\begin{center}
\renewcommand\arraystretch{1.1}
\begin{tabular}{|c|c|c||c|c|}
\hline
 & \multicolumn{2}{|c|}{ biarc PH curve} & \multicolumn{2}{|c|}{single PH curve}\\
\hline
$h$ & $E_{{\rm err}}$ & Decay exp.  & $E_{{\rm err}}$  & Decay exp.\\
\hline
$1$ & $5.23963\cdot 10^{-6}$ & /    &  $1.02470\cdot 10^{-5}$ & /\\
$ \frac{1}{2}$ & $1.48263 \cdot 10^{-7}$ & $5.14323$    &  $1.50891\cdot 10^{-7}$ &  $6.08555$ \\
$ \frac{1}{4}$ & $ 4.45464\cdot 10^{-9}$ & $5.05671$  &  $2.29413 \cdot 10^{-9}$ & $6.03941$ \\
$ \frac{1}{8}$ & $ 1.36921 \cdot 10^{-10}$ & $5.02389$  &  $3.53800\cdot 10^{-11}$ & $6.01887$\\
$ \frac{1}{16} $ & $ 4.24693 \cdot 10^{-12}$ & $5.01078$  & $5.49288 \cdot 10^{-13}$  & $6.00923$ \\
$ \frac{1}{32}$ & $ 1.32248 \cdot 10^{-13}$ & $5.00510$  & $8.55554\cdot 10^{-15}$  & $6.00456$ \\
$ \frac{1}{64}$ & $4.12567\cdot 10^{-15}$ & $5.00248$  & $1.33470\cdot 10^{-16}$  & $6.00227$\\
$ \frac{1}{128}$ & $ 1.28733 \cdot 10^{-16}$ & $5.00217$  &  $2.08384\cdot 10^{-18}$ & $6.00113$\\
$ \frac{1}{256}$ & $ 4.01726\cdot 10^{-18}$ & $5.00203$  & $3.25473 \cdot 10^{-20}$  & $6.00056$\\
\hline
\end{tabular}
 \end{center}
\caption{Example~\ref{Sec5-example5}. Errors $E_{\rm err}$ of  biarc interpolants (left part) and single PH curve interpolants (right part) for the curve $\f$ (with $\omega=0.2$)
on the interval $[0,h]$, together with estimates of the decay exponent as $h$ decreases to zero.
 }
\label{Tab:example-5-5}
\end{table}
\begin{figure}[htb]
\centering
\includegraphics[width=0.35\textwidth]{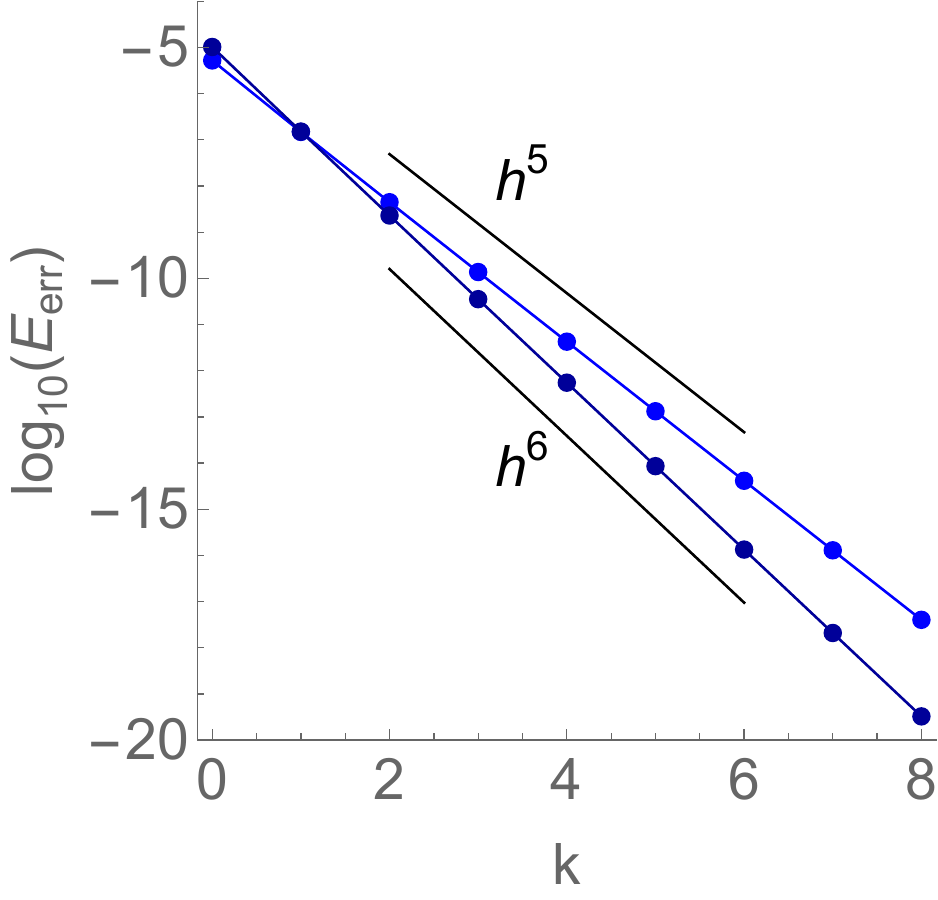}
\caption{
Example~\ref{Sec5-example5}. Errors $E_{\rm err}$ of biarc interpolants (light blue) and single PH curve interpolants (dark blue) for the curve $\f$ (with $\omega=0.2$)
on the interval $[0,2^{-k}]$ with respect to $k$ (in $\log_{10}$ scale).
}
\label{Fig:ex_5a}
\end{figure}

The presented interpolation scheme is completely local and yields $G^2$ PH spline curves when applied to the approximation of consecutive segments of the given curve $\f$. In particular, we choose
$[a,b]=[0,3\pi]$ and the uniform splitting of this interval by $s_j=\frac{3\pi j}{10}$ for $j=0,1,\dots,10$. Computing the PH biarcs that interpolate $\f$ on subintervals $[s_{j-1},s_{j}]$, $j=1,2,\dots,10$, we obtain the $G^2$ PH spline curve shown in Figure~\ref{Fig:ex_5b}, together with the porcupine plot of the curvature. The approximation error is computed as the maximum of errors  \eqref{eq:error} on each segment and equals  $2.09879\cdot10^{-5}$. To present it graphically
(see Figure~\ref{Fig:ex_5b})
we use an additional  linear reparameterization of the interval $[0,1]$ to $[j-1,j]$ for $j$-th segment, $j=1,2,\dots,10$.
In addition, we show  the difference between curvatures of the spline and the curve $\f$, which indicates that also the curvatures are well approximated.
\begin{figure}[htb]
\centering
\begin{minipage}{0.33\textwidth}
\includegraphics[width=1\textwidth]{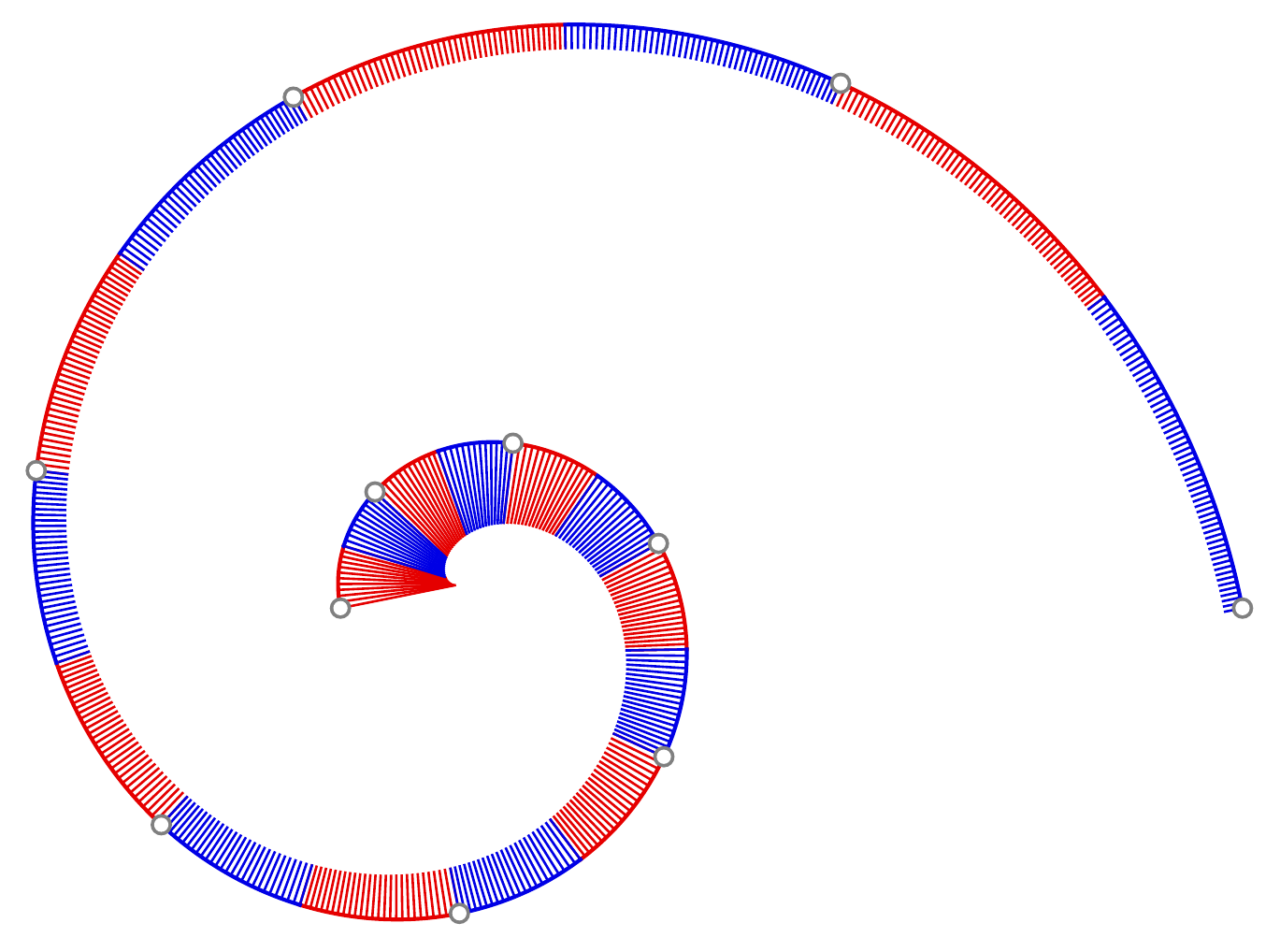}
\end{minipage}
\hskip.5cm
\begin{minipage}{0.28\textwidth}
\includegraphics[width=1\textwidth]{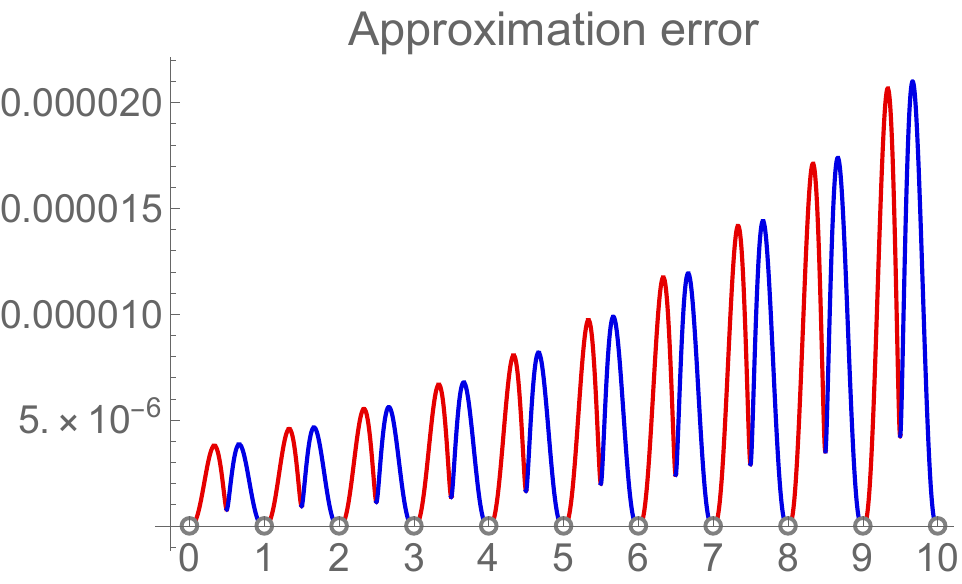}
\end{minipage}
\hskip.5cm
\begin{minipage}{0.28\textwidth}
\includegraphics[width=1\textwidth]{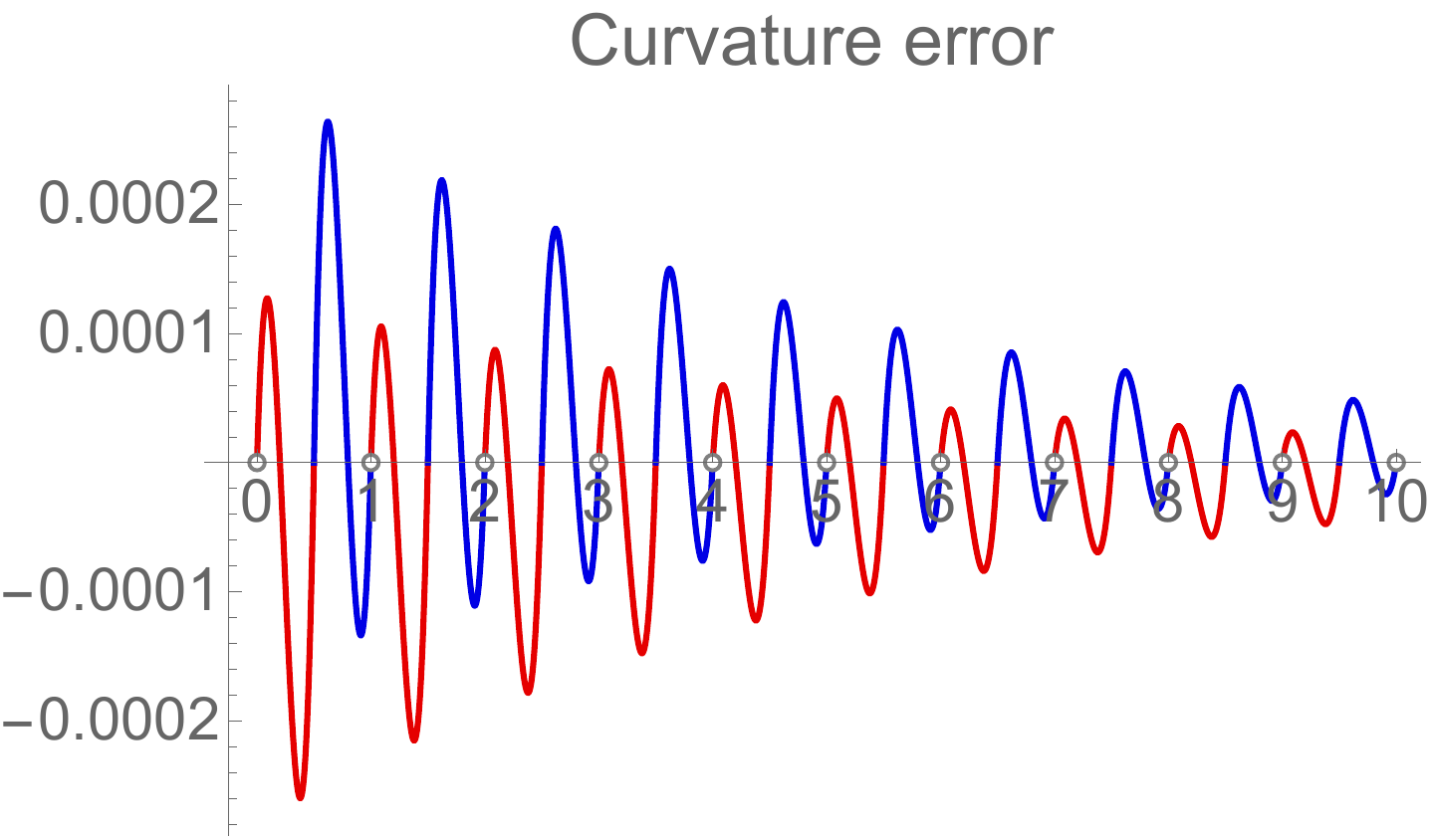}
\end{minipage}
\caption{
Example~\ref{Sec5-example5}.  The resulting $G^2$ PH spline interpolant equipped with the porcupine curvature plot (with no proportional factor), and graphs of approximation and curvature error.
}
\label{Fig:ex_5b}
\end{figure}
\end{exm}

As a final example let us consider the approximation of a circular arc, which is one of the most important objects in computer aided geometric design and there exist several results on its approximation with polynomials, see e.g.
\cite{VZ'19}, \cite{Vavpetic'20}, \cite{VZ'21}, and the references therein. Most of these papers propose methods which minimize the error and produce high order approximants. The main advantage of our PH biarc interpolant is that, in addition to a high order of approximation, the length of the circular arc is preserved while all the properties of PH curves can be applied, i.e. offset curves are rational, arc-length reparameterization is simple, etc.
\begin{exm}  \label{Sec5-example6}
One possible parameterization of the circular arc (of unit radius) follows from \eqref{eq:log-spiral} with $\omega=0$.
Choosing $a=0$, $b=\pi$, we obtain a semicircle and the corresponding PH biarc interpolants (the data are chosen by \eqref{data-from-f} for $h=\pi$) is shown in Figure~\ref{Fig:ex_6} (left), together with the porcupine curvature plot. It is obtained by $\beta_0=\beta_1=0$, $\lambda=1$
and $\alpha_0 = \alpha_1 = 1.77441$,  and approximates the semicircle with the error $E_{\rm err} = 6.38889\cdot10^{-4}$.
\begin{figure}[htb]
\centering
\begin{minipage}{0.35\textwidth}
\includegraphics[width=1\textwidth]{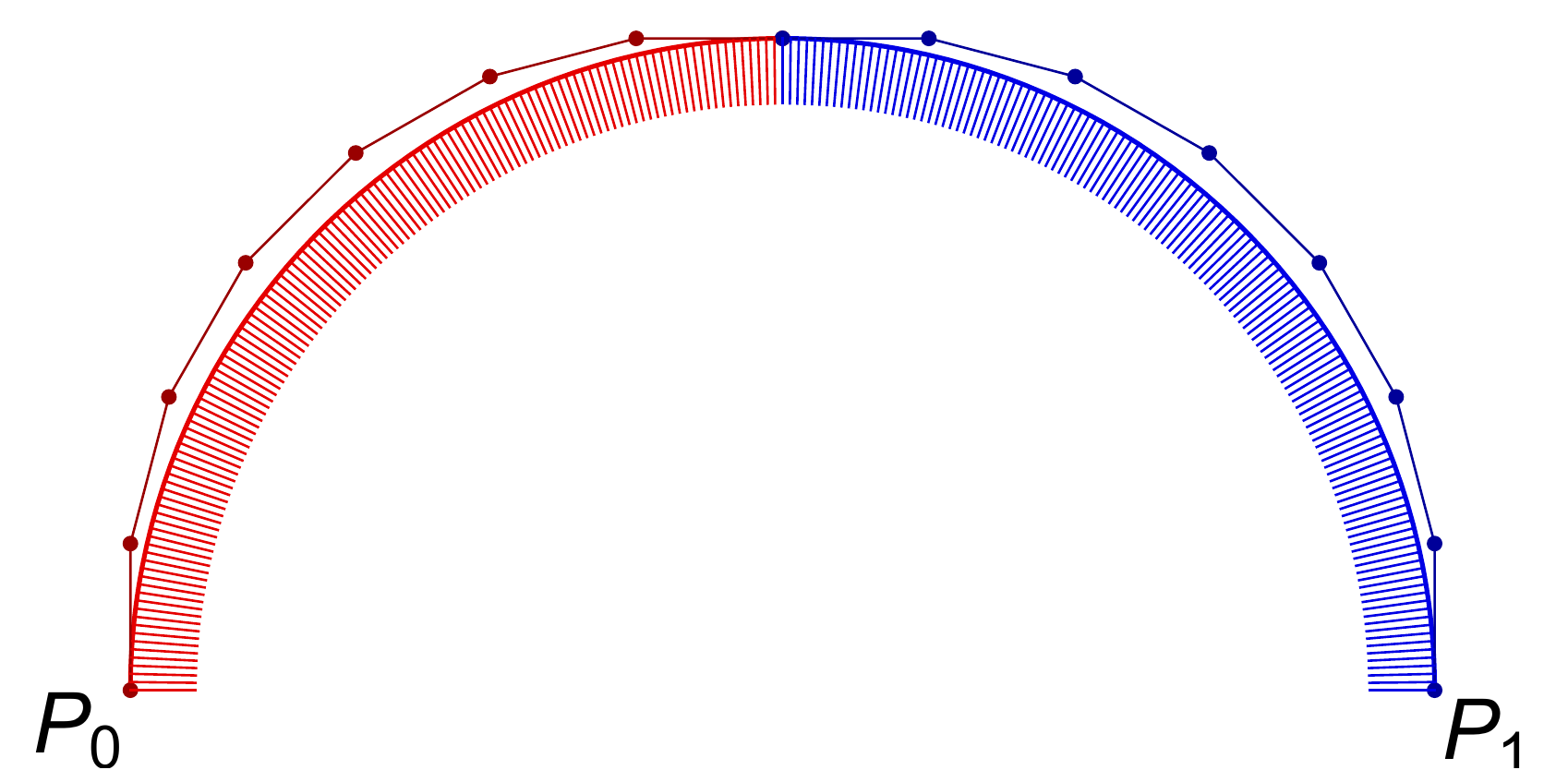}
\end{minipage}
\hskip.5cm
\begin{minipage}{0.28\textwidth}
\includegraphics[width=1\textwidth]{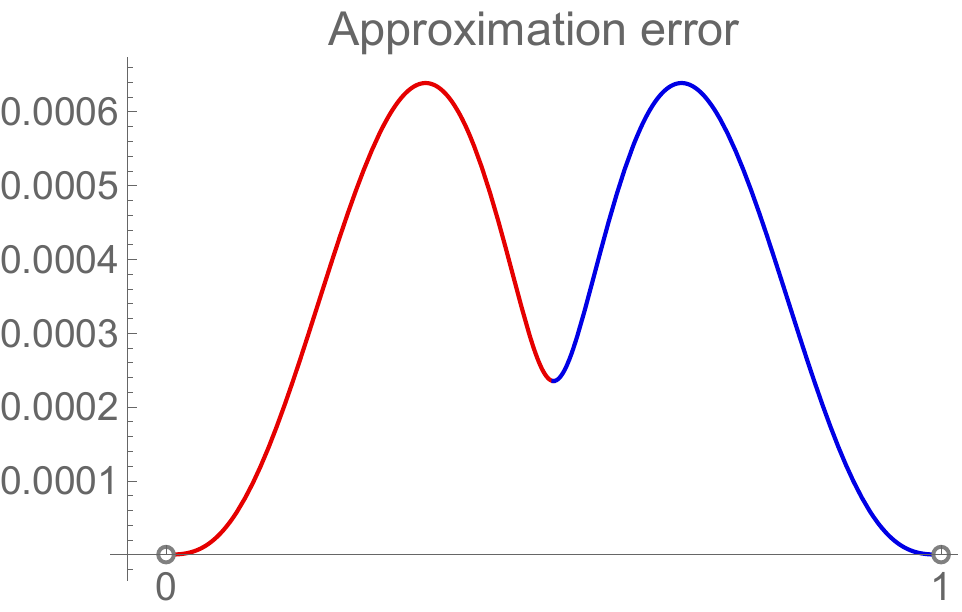}
\end{minipage}
\hskip.5cm
\begin{minipage}{0.28\textwidth}
\includegraphics[width=1\textwidth]{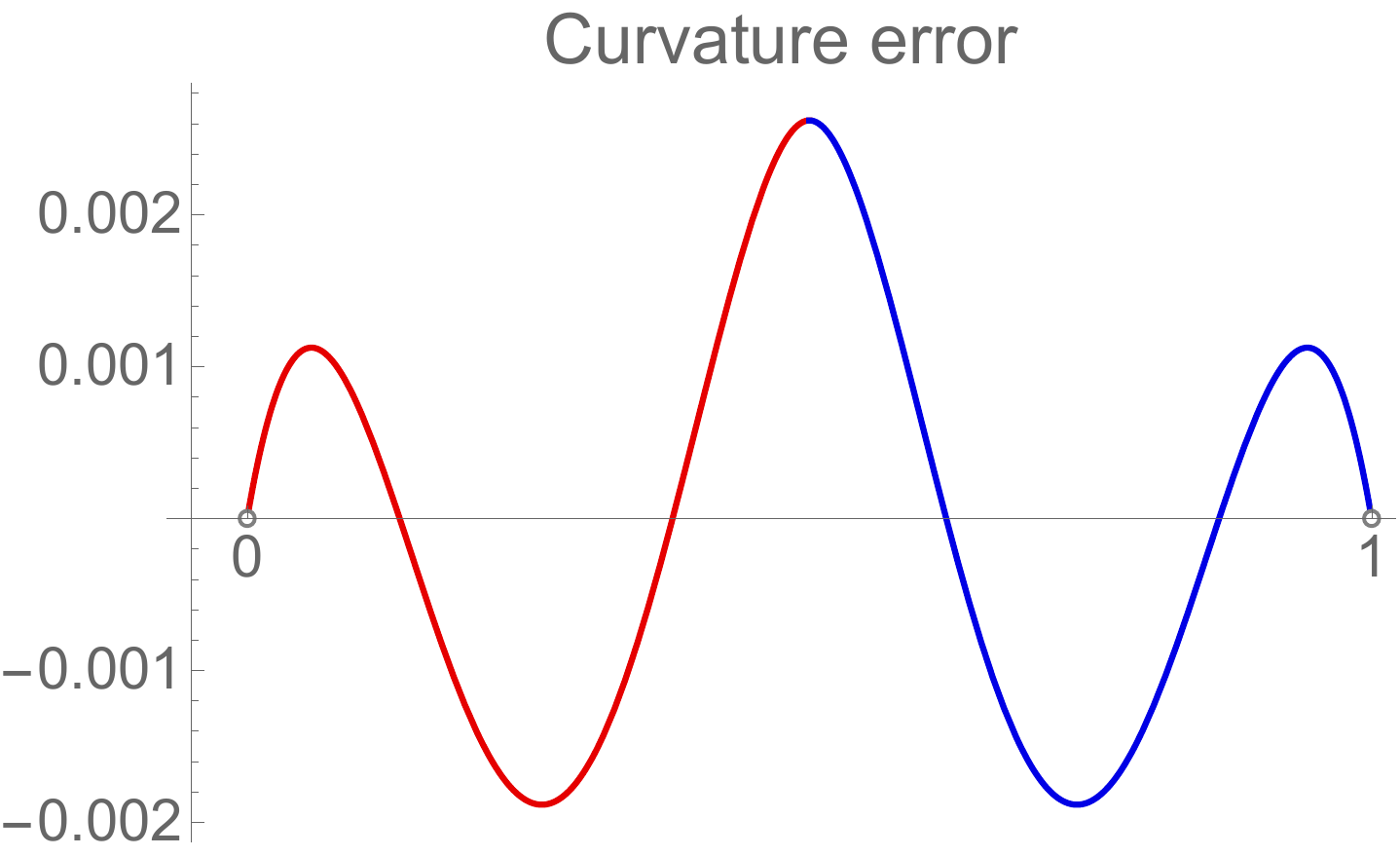}
\end{minipage}
\caption{
Example~\ref{Sec5-example6}.  The PH biarc that interpolates the semicircle, equipped with the porcupine curvature plot (with proportional factor equal to $0.1$), and graphs of approximation and curvature error.
}
\label{Fig:ex_6}
\end{figure}

Clearly, we obtain the $G^2$ spline interpolant of the whole circle by joining $N$ rotated versions of the PH biarc that approximates $\f$ over $[0, \phi]$ for $\phi=\frac{2\pi}{N}$. The approximation errors for increasing number of segments $N=2^k$, $k=1,2,\dots,9$, i.e. $\phi=2^{1-k} \pi$, are shown in Table~\ref{table-approxOrderCircle} and illustrated in Figure~\ref{fig-approxOrderCircle}, which indicates that the approximation order is $5$. \\
\begin{minipage}{0.45\textwidth}
\renewcommand\arraystretch{1.1}
\begin{tabular}{|c|c|c|c|}
\hline
$N$ & $\phi$ & $E_{{\rm err}}$ & {\rm Decay exp.}  \\
\hline
$2 $ & $\pi$ & $6.38885\cdot 10^{-4}$ & /    \\
$4 $ & $ \frac{\pi}{2}$ & $1.81754 \cdot 10^{-5}$ & $5.1355$    \\
$8 $ & $ \frac{\pi}{4}$ & $ 5.53785\cdot 10^{-7}$ & $5.03652$  \\
$16 $ & $ \frac{\pi}{8}$ &$1.71943 \cdot 10^{-8}$ & $5.00932$ \\
$32 $ & $ \frac{\pi}{16}$ & $ 5.36452 \cdot 10^{-10}$ & $5.00234$   \\
$64 $ & $ \frac{\pi}{32}$ & $ 1.67573\cdot 10^{-11}$ & $5.00059$  \\
$128 $ & $ \frac{\pi}{64}$ & $5.23612\cdot 10^{-13}$ & $5.00015$ \\
$256 $ & $ \frac{\pi}{128}$ & $  1.63625 \cdot 10^{-14}$ & $ 5.00004$  \\
$512$ & $ \frac{\pi}{256}$ & $ 5.11313\cdot 10^{-16}$ & $5.00004$  \\
\hline
\end{tabular}
\captionof{table}{Errors when interpolating the whole circle by the $G^2$ PH spline, composed of $N=2^k$ biarc segments. }
\label{table-approxOrderCircle}
\end{minipage}\hskip.5cm
\begin{minipage}{0.45\textwidth}
\includegraphics[width=0.8\textwidth]{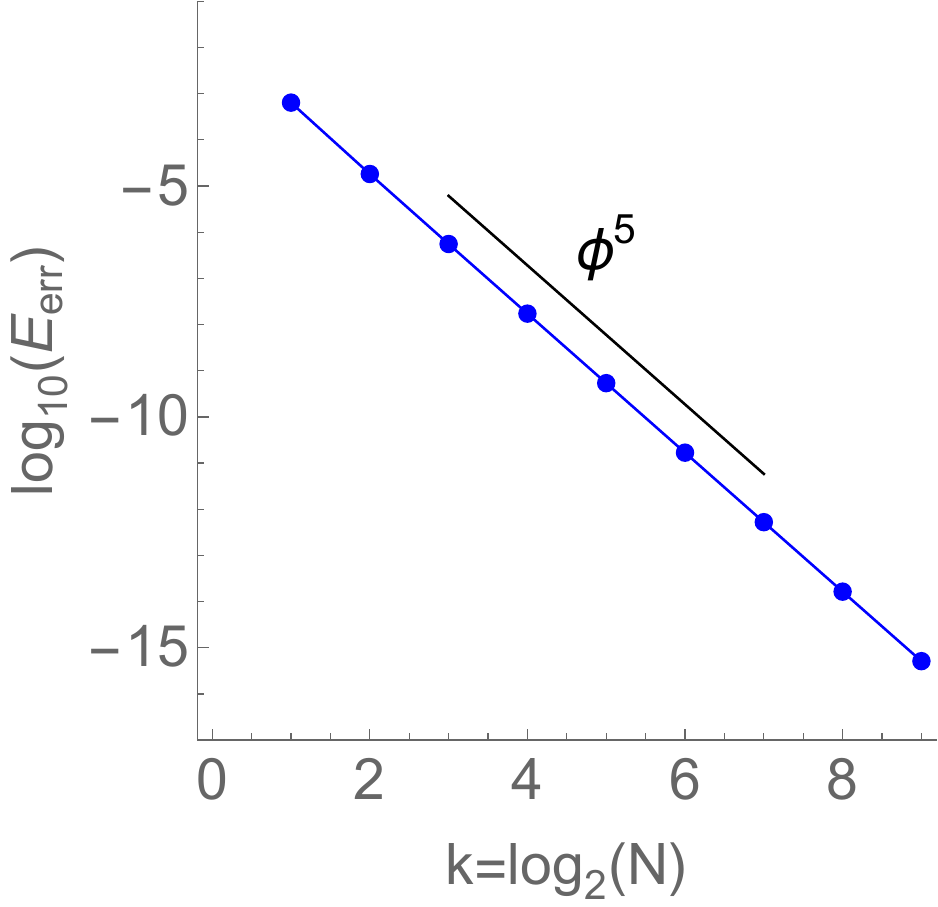}
\captionof{figure}{Graphical interpretation of the decay exponent. }
\label{fig-approxOrderCircle}
\end{minipage}
\end{exm}

\section{Conclusions}\label{sec:conclusion}
Most methods to construct curves rely on the interpolation of discrete data, such as points, tangents or curvatures. If, in addition, prescribed arc length is imposed, in general the use of iterative approximate methods is necessary.
Pythagorean--hodograph (PH) curves are polynomial curves with  the distinctive property of possessing arc lengths exactly determined by simple algebraic expressions in their coefficients. Hence the problem of constructing  $G^2$ planar curves, that interpolate points, tangent directions and curvatures, and in addition have prescribed arc-length, can be exactly addressed.
In this paper such problem is investigated considering PH curves of degree $7$ and it is shown that it reduces to three non-linear equations with one free shape parameter. However there exist data for which no interpolants can be found. A way to overcome this drawback is to consider {\it biarcs} keeping the degree to $7$. In this way
the solution of the $G^2$ continuity equations can be derived in a closed form,
depending on four free parameters. By fixing two of them to zero, it is proven that the length constraint can be satisfied for any data.  Beside of being easy to implement and simple to use in practice, the proposed method  can be directly applied to a (local) construction of $G^2$ continuous interpolating splines, as shown in  the final examples.

{
As a matter of further research, we aim to extend the proposed approach to the spatial case. We believe the extension is possible but not so straightforward.
Another interesting issue would be to consider the described $G^2$ interpolation problem using quintic PH biarcs. Counting the number of degrees of freedom, this could be possible, but the theoretical analysis of the existence of the solution is expected to be much more complicated.}

\section*{Acknowledgements}
 Research on this paper was supported in part by the program P1-0288 and the grant J1-7256 from ARRS, Republic of Slovenia, by the MIUR Excellence Department Project, awarded to the Department of
Mathematics, University of Rome ``Tor Vergata'' (CUP E83C18000100006F),
and by INdAM-GNCS, Gruppo Nazionale per il Calcolo Scientifico which F. Pelosi and M.L. Sampoli are members of.

\bibliography{PH-refs}

\end{document}